\documentclass{article}

\usepackage{delarray,verbatim,enumerate,a4wide}
\usepackage{amsmath,amsthm,amstext,amsbsy,amssymb,amsfonts,amscd}
\usepackage[all]{xy}
\usepackage[french]{babel}
\usepackage[T1]{fontenc}
\usepackage[utf8]{inputenc}

\usepackage{scalerel}

\usepackage{upgreek,bbm}

\usepackage[small]{titlesec}

\usepackage{color}
\definecolor{vert}{rgb}{0.1,0.4,0.2}
\usepackage[colorlinks=true,linkcolor=blue,citecolor=vert]{hyperref}

\usepackage{calligra}
\DeclareFontShape{T1}{calligra}{m}{n}{<->s*[0.95]callig15}{}
\DeclareMathAlphabet{\mathscr}{T1}{calligra}{m}{n}

\newtheorem{Th}{Théorème}[]
\newtheorem{Lem}[Th]{Lemme}
\newtheorem{Prop}[Th]{Proposition}
\newtheorem{Cor}[Th]{Corollaire}

\newtheorem{Sco}[Th]{Scolie}
\newtheorem{Conv} [Th]{Convention}
\newtheorem{Def} [Th]{Définition}

\newtheorem*{Th*}{Théorème}
\newtheorem*{Cor*}{Corollaire}
\newtheorem*{Def*}{Définition}

\def\Preuve{\noindent {\it Preuve.~}}

\def\Remarque{\smallskip\noindent {\it Remarque.~}}

\def\Nota{\smallskip\noindent {\it Nota.~}}

		\def\QQ{\mathbb Q}	
\def\NN{\mathbb N}	\def\ZZ{\mathbb Z}	\def\FF{\mathbb F}	
\def\F2{\mathbb{F}_2}	\def\Z2{\mathbb{Z}_2}		
\def\Zl{{\mathbb{Z}_\ell}} 	\def\Ql{{\mathbb{Q}_\ell}}	\def\Tl{\mathbb{T}_\ell}	\def\Fl{{\mathbb{F}_\ell}}

 				\def\U{\mathcal  U}	\def\F{\mathcal  F}	
  		\def\C{\mathcal  C}			\def\X{\mathcal  X}
 	  	\def\Cl{\mathcal  C\!\ell}	
\def\E{\mathcal  E}		\def\T{\mathcal  T}					\def\O{{\mathcal O}}

		\def\mmu{\boldsymbol{\upmu}}

		\def\p{{\mathfrak p}}				\def\a{{\mathfrak a}}		
		\def\l{{\mathfrak l}}						
		\def\h{{\mathfrak h}}		

			\def\Reg{{\R\!e\!g}}

\def\wi{\widetilde}				
			\def\dim{\operatorname{dim}}		\def\codim{\operatorname{codim}}
	\def\deg{\operatorname{deg}}				\def\In{\operatorname{In}}
\def\Gal{\operatorname{Gal}}					
				
\def\Reg{\operatorname{R\acute{e}g}}	\def\Aug{\operatorname{Aug}}		\def\Ram{\operatorname{Ram}}

\newcommand\scale[2]{\vstretch{#1}{\hstretch{#1}{#2}}}

\newcommand\si[1]{\scale{.7}{#1}}	
\newcommand\ph{{\phantom{*}}}

\def\%{{\scale{.8}{\infty}}}		\newcommand\et{{\scale{.8}{\rm \acute et}}}	\newcommand\amb{{\scale{.8}{\rm amb}}}

\makeatletter
\newcommand*\wt[2][0.2ex]{%
        \begingroup
        \mathchoice{\wt@helper{#1}{#2}{\displaystyle}{\textfont}}
                   {\wt@helper{#1}{#2}{\textstyle}{\textfont}}
                   {\wt@helper{#1}{#2}{\scriptstyle}{\scriptfont}}
                   {\wt@helper{#1}{#2}{\scriptscriptstyle}{\scriptscriptfont}}%
        \endgroup
        #2%
}
\newcommand*\wt@helper[4]{%
        \def\currentfont{\the#41}%
        \def\currentskewchar{\char\the\skewchar\currentfont}%
        \setbox\tw@\hbox{\currentfont$#2$\currentskewchar}%
        \dimen@ii\wd\tw@
        \setbox\tw@\hbox{\currentfont$#2${}\currentskewchar}%
        \advance\dimen@ii-\wd\tw@
        \rlap{\raisebox{-#1}{$\m@th#3\kern\dimen@ii\widetilde{\phantom{#2}}$}}%
}
\makeatother

			\def\wC{\wt[0.1ex]{\mathcal C}}

\begin{document}


\title{\LARGE \bf Genre des corps surcirculaires\footnote{Pub. Math. Besançon \no 2 (1986), 1--40.}}

\author{Jean-François {\sc Jaulent}}

\date{}
\maketitle\bigskip\bigskip

{\small
\noindent{\bf Résumé.}
Nous montrons que les formules de translation du genre à la Riemann-Hurwitz obtenues par Kuz'min, Kida, Iwasawa, Wingberg et alii pour l'invariant lambda attaché à certains modules d'Iwasawa d'une $\Zl$-extension cyclotomique de corps de nombres sont essentiellement équivalentes. Plus précisément, nous montrons que toutes ces formules, y compris celles énoncées à la Chevalley-Weil en termes de représentations, résultent identiquement pour des raisons purement algébriques du calcul arithmétique d'un quotient de Herbrand convenable qu'il suffit de mener dans le cas cyclique de degré premier.\medskip

\noindent {\bf Abstract.} We show that Riemann-Hurwitz-style translation formulas obtained by Kuz'min, Kida, Iwasawa, Wingberg et alii for the lambda invariant attached to certain Iwasawa moduli in cyclotomic $\Zl$-extension of number fields are essentially equivalent. More precisely, we prove that all these formulas, including those stated in terms of representations, result identically for purely algebraic reasons from the arithmetic computation of a suitable Herbrand quotient which it suffices to carry out in the cyclic case of prime degree $\ell$.\medskip

\noindent{\bf Avertissement} Le texte qui suit est la mise au format \LaTeX \,de l'article dactylographié original paru aux Publications Mathématiques de Besançon en 1986. Il n'en diffère que par la correction de diverses coquilles, par l'harmonisation de quelques notations avec celles des articles ultérieurs (notamment l'inversion de position de $S$ et $T$) ainsi que par la numérotation des théorèmes.

}


\tableofcontents

\section*{Introduction}
\addcontentsline{toc}{section}{Introduction}

On sait, depuis les travaux essentiels d'Iwasawa sur les corps cyclotomiques, que l'invariant $\lambda$ associé aux $\ell$-groupes de classes imaginaires dans une tour cyclotomique de corps à conjugaison complexe est, sous réserve de nullité de l'invariant $\mu$, un bon analogue pour la théorie des nombres de la notion classique de genre attachée aux corps de fonctions d'une variable\footnote{Une extension de corps imaginaires est dite admettre une conjugaison complexe lorsqu'elle provient d'une extension de même degré de leurs sous-corps réels}.\smallskip

Un exemple particulièrement frappant de cette analogie est la formule, démontrée par Kida \cite{Ki1}, qui relie les invariants $\lambda$ associés aux $\ell$-groupes de classes imaginaires dans une $\ell$-extension finie arbitraire de corps à conjugaison complexe et joue dans la théorie d'Iwasawa le rôle de la formule de Riemann-Hurwitz dans celle des corps de fonctions d'une variable. De fait, une formule semblable avait été établie peu avant par Kuz'min \cite{Ku} pour les invariants $\lambda$ associés aux $\ell$-groupes de $\ell$-classes imaginaires de ces corps (i.e. aux quotients des $\ell$-groupes de classes imaginaires par leurs sous-groupes respectifs engendrés par les classes des idéaux premiers au-dessus de $\ell$)\footnote{Le résultat de Kuz'min semble, aujourd'hui encore, curieusement méconnu. Il est vrai que la traduction anglaise de son long article aux Izvestija est légèrement postérieure au travail de Kida.}.
Ultérieurement Wingberg \cite{Wi0,Wi1} a montré qu'une formule identique vaut pour les invariants $\lambda$ associés aux $\ell$-groupes de classes infinitésimales\footnote{En fait, Wingberg travaille sur des groupes de Galois. L'interprétation infinitésimale provient de \cite{J9}.}.
Plus récemment des démonstration simplifiées de ces résultats ont été données par D'Mello et Madan \cite{DM} dans le cas étudié par Kida; par Gold et Madan \cite{GM1} dans celui considéré par Kuz'min.\smallskip

Auparavant, Iwasawa, remarquant que la formule de Kida se généralise en une identité entre représentations en tout point analogue à la formule de Chevalley-Weil pour les corps de fonctions (cf. \cite{CW}), en déduisit une preuve cohomologique du résultat de Kida \cite{Iw2}, que Gold et Madan \cite{GM2} ont récemment simplifiée.\smallskip
 
À ces démonstrations algébriques il convient d'ajouter, pour être complet, celles obtenues par voie analytique, qui utilisent la correspondance, étable par Iwasawa, entre invariants $\lambda$ et $\mu$, et fonctions $L$ $\ell$-adiques.La plus ancienne dans cette perspective est incontestablement la formule obtenue par Gras \cite{Gr0} pour les corps abéliens à l'aide des fonctions $L_\ell$ de Kubota-Leopoldt. La preuve générale de Sinnott \cite{Si} s'appuie sur la notion de pseudo-mesure $\ell$-adique introduite par Serre. Tout récemment Gold et Madan ont montré que les méthodes de Sinnott s'appliquent encore dans un cadre métacyclique où le degré des extensions n'est plus une puissance de $\ell$. Les relations entre fonctions $L$ qu'ils utilisent peuvent être regardées comme la version $\ell$-adique des identités complexes étudiées par Brauer et Walter (cf. \cite{J1}).\medskip

Cela dit, l'objet de cet article est triple:\smallskip

$(i)$ En premier lieu, nous montrons que les multiples formules évoquées plus haut ne sont que les diverses formulations d'une seule et même identité arithmétique. L'argument essentiel est ici la correspondance, exposée dans \cite{J15} et établie dans \cite{J14}, entre les invariants $\lambda$ associés aux différents groupes étudiés, qui affirme que tous les paramètres considérés se ramènent, sous les conjectures de Leopoldt et de Gross-Kuz'min, à l'un quelconque d'entre eux au moyen de formules standard.\smallskip

$(ii)$ Nous établissons ensuite que les deux expressions possibles de cette identité (en termes numériques ou bien en termes de représentations, c'est-à-dire
de caractères) sont rigoureusement équivalentes. Ce point résulte d'un argument classique de Herbrand, qui ramène l'étude d'une représentation galoisienne à la détermination des restrictions de son caractère aux sous-groupes cycliques de son groupe de Galois.\smallskip

$(iii)$ Enfin, nous prouvons que, par une heureuse coïncidence du vocabulaire, l'identité en question relève intégralement de la théorie habituelle des genres pour les corps de nombres. Nous en donnons deux démonstrations algébriques : la première, écrite en termes de groupes de classes, s'obtient simplement par passage à la limite inductive à partir des suites exactes classiques de la théorie des genres; la seconde, écrite en termes de groupes de Galois, s'obtient par passage à la limite projective. Nous renvoyons à \cite{Gr2} pour une théorie analytique des genres basée sur la notion de pseudo-mesure $p$-adique.

\newpage
\section{Description des paramètres lambda attachés à un corps surcirculaire}

\subsection{Présentation des corps surcirculaires}

Le nombre premier $\ell$ étant fixé une fois pour toutes, nous désignons par $\QQ_\%$ la $\Zl$-extension cyclotomique du corps des rationnels.
Les invariants que nous nous proposons d'étudier sont attachés à une catégorie bien particulière de corps, stable par extension finie:

\begin{Def}\label{D1}
Nous appelons corps surcirculaire (relativement au nombre premier $\ell$) toute extension algébrique finie de la $\Zl$-extension  cyclotomique $\QQ_\%$ du corps des rationnels.
\end{Def}

\Remarque Iwasawa appelle $\Zl$-corps les extensions obtenues comme $\Zl$-extensions d'un corps de nombres (de degré fini). Cette notion ne nous parait pas exactement adaptée au problème que nous considérons, puisque les résultats que nous avons en vue supposent explicitement que la $\Zl$-extension étudiée soit la cyclotomique.
\medskip

Tout corps surcirculaire $K_\%$ peut être regardé d'une infinité de façons comme la $\Zl$-extension cyclotomique d'un corps de nombres $K$.
Par exemple, si $x$ est un élément primitif de $K_\%$ sur $\QQ_\%$, il est immédiat que $K=\QQ[x]$ convient.
Quoique non canonique, cette description va nous être utile pour appliquer les théorèmes de structure de la théorie d'Iwasawa.\smallskip

Notons auparavant qu'en vertu de la loi de décomposition des idéaux premiers dans une tour cyclotomique, le fait que tout corps surcirculaire s'écrive $K\QQ_\%$ pour un corps de nombres convenable $K$, montre que les places ultramétriques de $\QQ$ sont finiment décomposées dans $K_\%/\QQ$.
Plus précisément, nous avons le résultat bien connu suivant:

\begin{Lem}\label{L1}
Soient $K_\%$ un corps surcirculaire et $K$ un corps de nombres vérifiant $K_\%=K\QQ_\%$.
\begin{itemize}
\item[(i)] Les places de $K$ au-dessus de $\ell$ sont presque totalement ramifiées dans $K_\%/K$; c'est-à-dire qu'il existe une sous-extension $K_{n_{\si{0}}}$ de $K_\%$, de degré fini $\ell^{n_{\si{0}}}$ sur $K$, telle que les places au-dessus de $\ell$ soient totalement ramifiées dans $K_\%/K_{n_{\si{0}}}$.
\item[(ii)] Les autres places ultramétriques sont presque totalement inertes dans $K_\%/K$; c'est-à-dire qu'il existe, pour chaque premier $\p$ de $K$, une sous-extension $K_{n_{\si{\p}}}$ de $K_\%$, de degré fini $\ell^{n_{\si{\p}}}$ sur $K$, telle que les places au-dessus de $\p$ soient totalement inertes dans $K_\%/K_{n_{\si{\p}}}$.
\item[(iii)] Les places archimédiennes sont complètement décomposées dans $K_\%/K$.
\end{itemize}
\end{Lem}

Une conséquence élémentaire mais peu citée de l'assertion $(ii)$ concerne les corps résiduels: si $\p$ est une place de $K$ au-dessus d'un premier $p$ différent de $\ell$, le corps résiduel $F_{\p_\%}=\O_\%/\p_\%$ associé à l'une quelconque $\p_\%$ des places de $K_\%$ au-dessus de $\p$ est la $\Zl$-extension (cyclotomique) de $F_\p$. En particulier, dans une $\ell$-extension de corps surcirculaires $L_\%/K_\%$, il ne peut y avoir d'inertie aux places ultramétriques étrangères à $\ell$.\smallskip

Enfin, il résulte du Lemme que le groupe des diviseurs d'un corps surcirculaire peut se décrire comme suit:
\begin{itemize}
\item[(i)] Les diviseurs construits sur les places au-dessus de $\ell$ forment un $\QQ$-espace vectoriel de dimension finie:

\centerline{$D_{K_\%}(\ell)=\bigoplus_{\p_\%\mid p}\;\p_\%^\QQ$}\smallskip

\item[(ii)] Les diviseurs construits sur les places au-dessus d'un premier $p\ne\ell$ forment un $\ZZ$-module libre de dimension finie:

\centerline{$D_{K_\%}(p)=\bigoplus_{\p_\%\mid\ell}\;\p_\%^\ZZ$}\smallskip

\item[(iii)] Les diviseurs construits sur les places archimédiennes forment un $\FF_2$-espace vectoriel qui est soit nul (si $K$ n'admet aucun plongement réel), soit infini dénombrable (si $K$ admet un et donc une infinité de plongements réels):

\centerline{$D_{K_\%}(\infty)=\bigoplus_{\p_\%\mid \infty}\;\p_\%^{\ZZ/2\ZZ}$.}
\end{itemize}\medskip

Et il vient alors: \qquad $D_{K_\%}=D_{K_\%}(\ell)\bigoplus \big(\bigoplus_{p\ne\ell}D_{K_\%}(p)\big)\bigoplus D_{K_\%}(\infty)$.

\subsection{Les théorèmes de structure de la théorie d'Iwasawa}

Reprenons la description des corps surcirculaires donnée plus haut.
Partons d'un corps de nombres $K$ (de degré fini sur $\QQ$; et considérons la $\Zl$-extension cyclotomique $K_\%$ de $K$.\par

Notons $\Gamma=\Gal(K_\%/K)$ le groupe de Galois correspondant, identifié à $\Zl$ par le choix d'un générateur topologique $\gamma$; et $\Lambda=\Zl[[\gamma-1]]$ l'algèbre d'Iwasawa associée. Écrivons de même $K_n$ l'unique sous-corps de $K_\%$ cyclique de degré $\ell^n$ sur $K$ et $\Gamma_{\!n}=\Gal(K_\%/K_n)$ le sous-groupe de $\Gamma$ associé à $K_n$ par la théorie de Galois.\smallskip

Supposons choisis deux ensembles finis disjoints $S$ et $T$ de places de $\QQ$ et considérons les (pro)-$\ell$-groupes $\,\Cl_T^S(K_n)$ de $T$-classes $S$-infinitésimales des corps $K_n$\footnote{La définition générale précise des groupes $\,\Cl^S_T$ est sans importance ici, car nous ne serons intéressés que par des spécialisations très particulières de $S$ et $T$. Pour la définition et leurs principales propriétés, cf. \cite{J18}, Ch. II, \S 1.}.
Disons simplement ici que $\,\Cl^S_T(K_n)$ est le (pro)-$\ell$-groupe de classes de diviseurs qui s'identifie canoniquement, via la théorie du corps de classes, au groupe de Galois $G^S_T(K_n)$ de la (pro)-$\ell$-extension abélienne maximale de $K_n$ qui est $S$-ramifiée et $T$ décomposée (c'est-à-dire non ramifiée aux places étrangères à $S$, complètement décomposée aux places au-dessus de $T$, et maximale sous ces deux conditions).\smallskip

Au sommet de la tour cependant, cet isomorphisme n'a plus lieu et il importe de distinguer entre groupes de classes et groupes de Galois. C'est pourquoi nous posons:

\begin{Def}\label{Not}
Étant donné un corps surcirculaire $K_\%=\bigcup_{n\in\NN}K_n$ et deux ensembles finis disjoints $S$ et $T$ de places de $\QQ$, nous écrivons:

$(i)$ $\,\C^S_T(K_\%)=\varprojlim \,\Cl^S_T(K_n)$  la limite du système projectif des $\,\Cl^S_T(K_n)$ pour la norme arithmétique.
Le groupe $\,\C^S_T(K_\%)$ est un $\Lambda$-module compact, canoniquement isomorphe au groupe de Galois $G^S_T(K_\%)$ de la pro-$\ell$-extension abélienne $S$-ramifiée $T$-décomposée maximale de $K_\%$.\smallskip

$(ii)$ $\,\Cl^S_T(K_\%)=\varinjlim \,\Cl^S_T(K_n)$ la limite du système inductif des $\,\Cl^S_T(K_n)$ pour l'extension des diviseurs. Et nous disons que $\,\Cl^S_T(K_\%)$ est le $\ell$-groupe des $T$-classes $S$-infinitésimales de $K_\%$.
\end{Def}

Les théorèmes de structure de la théorie d'Iwasawa affirment alors l'existence de trois invariants $\rho^S_T$, $\mu^S_T$ et $\lambda^S_T$ qui mesurent la croissance avec $n$ des quotients d'exposant $\ell^n$ respectifs des groupes $\,\Cl^S_T(K_n)$ et d'une constante $\nu^S_T$ (éventuellement négative) telle qu'on ait pour $n>>0$:\smallskip

\centerline{$|^{\ell^{\si{n}}\!}\Cl^S_T(K_n)|=\ell^{x_{\si{n}}}$, avec $x_n=\rho^S_T n\ell^n+\mu^S_T\ell^n+ \lambda^S_T n+ \nu^S_T$.}\smallskip

\noindent Les paramètres $\rho^S_T$ et $\mu^S_T$, qui dépendent du choix de $K$, ne sont pas à proprement parler des invariants de $K_\%$ contrairement à l'invariant $\lambda^S_T$, qui est au cœur de cette étude. Précisons cela:

\begin{Prop}\label{P1}
Le groupe $\,\C^S_T(K_\%)$ est un $\Lambda$-module de torsion si et seulement si l'une ou l'autre des deux conditions suivantes est vérifiée:
\begin{itemize}
\item[(i)] ou bien l'ensemble $T$ ne contient pas la place $\ell$,
\item[(ii)] ou bien $T$ contient $\ell$ et $K$ est totalement réel.
\end{itemize}
Dans tous les autres cas, le paramètre $\rho^S_T$ est égal au nombre de places complexes du corps $K$.
\end{Prop}

\Preuve Ce résultat étant essentiellement connu, contentons-nous d'une esquisse de démonstration:

$(i)$ Si $T$ ne contient pas $\ell$, les groupes $\,\C^S_T(K_n)$ sont finis et leur limite projective est donc $\Lambda$-module de torsion par un argument classique de théorie d'Iwasawa.

$(ii)$ Si $T$ contient $\ell$, les groupes $\,\C^S_T(K_n)$ sont infinis et deux cas se présentent:
\begin{itemize}
\item Lorsque le corps $K$ n'est pas totalement réel, il possède au moins un conjugué imaginaire et $K_\%$ admet alors une infinité de $\Zl$-extensions indépendantes, lesquelles sont non-ramifiées en dehors de $\ell$ et $T$-décomposées (car il ne peut y avoir inertie aux places ultramétrique étrangères à $\ell$). Et le groupe $\,\C^S_T(K_\%)$ n'est donc pas un $\Lambda$-module de torsion.
\item Lorsque, en revanche, $K$ est totalement réel,$K_\%$ n'admet qu'un nombre fini de $\Zl$-extensions indépendantes $\ell$-ramifiées, en vertu de l'exactitude de la conjecture faible de Leopoldt dans les tours cyclotomiques (cf. e.g. \cite{J18}, Ch. IV, \S 2). Le groupe $\,\C^{\{\ell\}}_T(K_\%)$ est ainsi un $\Lambda$-moduls de torsion. Et il en va de même pour le groupe $\,\C^S_T(K_\%)$, puisque les places de $K_\%$ sont en nombre fini et que le sous-module de $\,\C^S_T(K_\%)$ engendré par les sous-groupes d'inertie associés  est donc un $\Zl$-module de type fini.
\end{itemize}

Dans les deux cas, le paramètre $\rho^S_T$ est celui du groupe $\,\C^{\{\ell\}}_\emptyset (K_\%)$ qui correspond à la pro-$\ell$-extension abélienne $\ell$-ramifiée maximale du corps $K_\%$. Et comme il est bien connu que celui-ci est égal au nombre de places complexes du corps $K$\footnote{L'égalité $\rho^{\{\ell\}}_\emptyset=c_K$ est, en fait, équivalente à l'exactitude de la conjecture faible de Leopoldt dans $K_\%/K$.}, la proposition est ainsi établie.

\begin{Prop}\label{P2}
Lorsque le corps $K$ contient les racines $2\ell$-ièmes de l'unité, les groupes $\,\C^S_T$ ont même paramètre $\mu$ quels que soient les ensembles finis $S$ et $T$ de places de $\QQ$. En particulier, l'égalité $\mu^S_T=0$ a lieu, en dehors de toute autre hypothèse sur $K_\%$, dès que le groupe de Galois $G^\emptyset_{\{\ell\}}(K_\%[\zeta_{2\ell}])$ associé à la pro-$\ell$-extension abélienne non-ramifiée $\ell$-décomposée maximale de corps surcirculaire $K_\%[\zeta_{2\ell}]$ est un $\Zl$-module de type fini\footnote{On ne connait pas de corps surcirculaire dans lequel cette condition est en défaut.}.
\end{Prop}

\Preuve
Les résultats de \cite{J14} montrent que les groupes $\,\C(K_\%)=C^{\,\emptyset}_{\,\emptyset}(K_\%)$, $\,C_{\{\ell\}}(K_\%)=C^{\,\emptyset}_{\{\ell\}}(K_\%)$ et $\,C^{\{\ell\}}(K_\%)=C_{\,\emptyset}^{\{\ell\}}(K_\%)$, qui correspondent respectivement à la pro-$\ell$-extension abéliennes non-ramifiée maximale, à celle qui est non-ramifiée et $\ell$-décomposée et à celle qui est $\ell$-ramifiée, ont même paramètre $\mu$ dès lors que le corps $K_\%$ contient le groupe $\mmu_{\ell^\%}$ des racines de l'unité d'ordre $\ell$-primaire\footnote{Pour une autre démonstration de ce résultat, voir \cite{J18}, Ch. IV, \S2.}.
La proposition en résulte immédiatement, puisque les places de $S$ ou $T$ qui sont étrangères à $\ell$ sont sans influence sur le paramètre mu, les une parce qu'il ne peut y avoir d'inertie dans une $\ell$-extension de $K_\%$ en une place étrangère à $\ell$, les autres parce le sous-groupe d'inertie d'une telle place dans une pro-$\ell$-extension de $K_\%$ est un $\Zl$-module de type fini\footnote{Cet argument serait en défaut dans un $\Zl$-corps non-surcirculaire, où les places de $S$ étrangères à $\ell$ seraient susceptibles d'être infiniment décomposées.}.\medskip

Iwasawa a conjecturé la nullité du paramètre $\mu$ attaché au groupe $\,\C(K_\%)$. D'après la Proposition \ref{P2}, cela revient à postuler que l'on a: $\mu_S^T(K_\%)=0$ pour tout corps surcirculaire $K_\%$ quels  que soient les ensembles finis $S$ et $T$ de places de $\QQ$. Il en et bien ainsi, par le théorème de Ferrero et Washington, dès que le corps $K_\%$ est abélien sur $\QQ$ (cf. e.g. \cite{Wa}, Th. 7.15). Mais ce n'est plus vrai en général pour les $\Zl$-corps non-surcirculaires (cf. \cite{Iw0}).

\subsection{Le problème de la capitulation}

Nous avons dit plus haut que les valeurs exactes des paramètres $\rho^S_T$ et $\mu^S_T$ attachés à un corps surcirculaire $K_\%$ dépendent du choix du corps de base $K$. Il n'en est pas de même, en revanche, des conditions  $\rho^S_T=0$ et $\mu^S_T=0$.\smallskip
\begin{itemize}
\item[(i)] La première $\rho^S_T=0$ signifie, en effet, que le $\Ql$-espace vectoriel $\Ql\otimes_\Zl \C^S_T(K_\%)$ construit sur le $\Zl$-module $\,\C^S_T(K_\%)$ est de dimension finie.
\item[(ii)] La seconde $\mu^S_T=0$ affirme, elle, que le sous-module de $\Zl$-torsion de $\,\C^S_T(K_\%)$ est fini.
\end{itemize}\smallskip

Lorsqu'elles sont vérifiées, le groupe $\,\C^S_T(K_\%)$ est ainsi un $\Zl$-module de type fini, somme directe comme tel d'un $\Zl$-module libre de dimension $\lambda_T^S$ et du sous-groupe de $\Zl$-torsion $\T^S_T(K_\%)$:\smallskip

\centerline{$\C^S_T(K_\%) \simeq \Zl^{\lambda_T^S}\oplus\T^S_T(K_\%)$.}\smallskip

Dans ce cas, les résultats de \cite{GJ} interprètent le sous-groupe $\T^S_T(K_\%)$ comme la limite projective des noyaux $\C\!ap^S_T(K_n)$ des applications naturelles  $j_n$ des $\ell$-groupes de $T$-classes $S$-infinitésimales $\,\Cl^S_T(K_n)$ des corps $K_n$ dans leur limite inductive $\,\Cl^S_T(K_\%)$, laquelle est dans les cas considérés ici soit un $\Zl$-module divisible de codimension $\lambda^S_T$, soit (sous la conjecture de Leopoldt) somme directe d'un tel module et d'un exemplaire de $\Zl$:\smallskip
\begin{itemize}
\item[(i)] Le premier cas $\,\Cl^S_T\simeq(\Ql/\Zl)^{\lambda^S_T}$ se produit lorsque l'ensemble $S$ ne contient pas $\ell$. Les groupes $\,\Cl^S_T(K_n)$ sont alors finis, et les arguments de \cite{GJ} donnent directement le résultat.
\item[(ii)] Le second cas $\,\Cl^S_T\simeq (\Ql/\Zl)^{\lambda^S_T} \oplus \Zl$ se produit lorsque $S$ contient $\ell$. L'hypothèse $\rho^S_T=0$ implique que les corps $K_n$ sont totalement réels, de sorte que les groupes $\,\Cl^S_T(K_n)$ s'écrivent, sous la conjecture de Leopoldt, comme somme directe de leurs sous-modules de torsion respectifs $\,\T^S_T(K_n)$ et d'un $\Zl$-module libre de dimension 1, isomorphe à $\Gal(K_\%/K_n)$.
Ce dernier module disparait par passage à la limite projective, mais non par passage à la limite inductive, de sorte que $\,\Cl^S_T(K_\%)$ est bien, comme annoncé, somme directe d'un $\Zl$-module libre de dimension 1 et de la limite inductive des $\,\T^S_T(K_n)$, qui est un $\Zl$-module divisible de codimension $\lambda_T^S$ par les mêmes arguments que plus haut.
\end{itemize}\smallskip

Les groupes $\C\!ap^S_T(K_n)$ mesurent ce qu'il est convenu d'appeler, depuis Hilbert, une capitulation; en l’occurrence la $T$-capitulation $S$-infinitésimale dans l'extension profinie $K_\%/K_n$. Avouons tout de suite qu'ils sont en général très mal connus et que c'est là un des obstacles essentiel au calcul de l'invariant $\lambda^S_T$. La seule information directe que nous possédons sur eux provient de la suite exacte des classes ambiges de Chevalley, en fait de sa généralisation aux pro-$\ell$-groupe de $T$-classes $S$-infinitésimales (cf. \cite{J18}, Ch. II, \S2), qui interprète $\C\!ap^S_T(K_n)$ comme un sous-groupe du premier groupe de cohomologie $H^1(\Gamma_{\!n},\E_T^S(K_\%))$ construit sur les $T$-unités $S$-infinitésimales dans l'extension profinie $K_\%/K_n$.
Faute de mieux, on est alors amené à restreindre l'étude des invariants $\lambda^S_T$ aux seuls cas où le groupe $H^1(\Gamma_{\!n},\E_T^S(K_\%))$ est trivial pour des raisons évidentes:
\begin{itemize}
\item[(i)] Lorsque $S$ contient $\ell$ (auquel cas les hypothèses faites implique que $K_\%$ est totalement réel), l'exactitude de la conjecture faible de Leopoldt montre que les groupes $\E_T^S(K_n)$ sont ultimement constants\footnote{La conjecture faible de Leopoldt revient à affirmer que les pro-$\ell$-groupes d'unités infinitésimales sont ultimement constants dans une $\Zl$-extension cyclotomique.(cf. e.g. \cite{J14}). L'extension de ce résultat aux groupes de $T$-unités $S$-infinitésimales est sans malice.}, ce qui entraîne banalement $H^1(\Gamma_{\!n},\E_T^S(K_\%)) = H^1(\Gamma_{\!n},\E_T^S(K_n))=1$ d'où  $\C\!ap^S_T(K_n)=1$, comme attendu.
\item[(ii)] Lorsque, en revanche, $S$ ne contient pas $\ell$, la cohomologie des groupes $\,\E_T^S$ n'est pas connue en général, ce qui amène à faire des hypothèses restrictive sur $K_\%$; dans la pratique, à se restreindre aux seules composantes imaginaires:
\end{itemize}

\begin{Def}\label{D2}
Nous disons qu'une extension algébrique de $\QQ$ admet une conjugaison complexe lorsque c'est une extension quadratique totalement imaginaire d'un sous-corps totalement réel. Un corps surcirculaire à conjugaison complexe est donc une extension quadratique totalement imaginaire d'un corps surcirculaire totalement réel.
\end{Def}

En présence d'une conjugaison complexe $\bar\tau$, il est naturel de dire qu'un module multiplicatif est réel lorsque $\bar\tau$ opère trivialement dessus; qu'il est imaginaire lorsqu'elle opère par passage à l'inverse. Si le nombre premier $\ell$ est impair, tout $\Zl$-module $M$ s'écrit ainsi comme somme directe:\smallskip

\centerline{$ M=M^{\si{\oplus}}\oplus M^{\si{\ominus}}=M^{(1+\bar\tau)/2}\oplus M^{(1-\bar\tau)/2}$}\smallskip

\noindent de ses composantes réelles et imaginaires. Cependant, si $\ell$ vaut 2, la même somme $M^{\si{\oplus}} + M^{\si{\ominus}}$ peut n'être ni directe, ni totale. Cela n'interdit pas, bien sûr, de continuer à parler de composantes réelles ou imaginaires,mais requiert à l'occasion quelques précautions. C'est pourquoi, nous supposons dans tout ce qui suit que $\ell$ est impair.\smallskip

Ce point acté, l'intérêt de se restreindre aux composantes imaginaires est que dans une tour cyclotomique e corps à conjugaison complexe, il n'est d'autres unités imaginaires que les racines de l'unité présentes dans la tour, lesquelles sont cohomologiquement triviales. Cela suffit à entraîner la nullité des groupes $H^1(\Gamma_{\!n},\E^S_T(K_n))^{\si{\ominus}}$ pour $n\gg 0$, et donc celle de la composante imaginaire $\C\!ap_T^S(K_\%)^{\si{\ominus}}$ du sous-groupe de torsion de $\C^S_T(K_\%)$.

\begin{Prop}\label{P3}
Supposons que $\ell$ soit impair, que le corps surcirculaire $K_\%$ admette une conjugaison complexe et que les invariants $\mu^S_T$ soient nuls.
\begin{itemize}
\item[(i)] Si $S$ contient la place $\ell$, la composante réelle du groupe $\,\C^S_T(K_\%)$ est un $\Zl$-module libre:\smallskip

\centerline{$\C^S_T(K_\%)^{\si{\oplus}} \simeq \Zl^{{\lambda^S_T}^{\si{\oplus}}}$.}\smallskip

Et si la conjecture de Leopoldt est vérifiée dans $K_\%$ (pour $\ell$), la composante réelle du groupe $\,\Cl^S_T(K_\%)$ est la somme directe d'un exemplaire de $\Zl$ et d'un $\Zl$-module divisible:\smallskip

\centerline{$\Cl^S_T(K_\%)^{\si{\oplus}} \simeq \Zl \oplus (\Ql/\Zl)^{{\lambda^S_T}^{\si{\oplus}}}$.}\smallskip 

\item[(ii)] Si $S$ ne contient pas $\ell$, la composante imaginaire de $\,\C^S_T(K_\%)$ est un $\Zl$-module libre:\smallskip

\centerline{$\C^S_T(K_\%)^{\si{\ominus}} \simeq \Zl^{{\lambda^S_T}^{\si{\ominus}}}$.}\smallskip

Et dans ce cas la composante imaginaire de $\,\Cl^S_T(K_\%)$ est un $\Zl$-module divisible:\smallskip

\centerline{$\Cl^S_T(K_\%)^{\si{\ominus}} \simeq (\Ql/\Zl)^{{\lambda^S_T}^{\si{\ominus}}}$.}
\end{itemize}
\end{Prop}

\subsection{Correspondance entre les divers paramètres lambda}

Pour voir que les paramètres lambda correspondant aux divers choix de $S$ et de $T$ se ramènent tous à l'un quelconque d'entre eux dans chacun des trois cas fondamentaux évoqués plus haut, nous pouvons nous appuyer sur les résultats de \cite{J18}:\smallskip

\begin{itemize}

\item[(i)] Le cas $S=\emptyset$, $T=\emptyset$ est celui considéré par Kida \cite{Ki1}: le groupe $\,\C(K_\%)=\C^{\,\emptyset}_{\,\emptyset}(K_\%)$ s'identifie au groupe de Galois de la pro-$\ell$-extension non-ramifiée maximale de $K_\%$; et le $\ell$-groupe de classes associé $\,\Cl(K_\%)=\Cl^{\,\emptyset}_{\,\emptyset}(K_\%)$ n'est autre que le $\ell$-groupe des classes de diviseurs au sens ordinaire du corps $K_\%$. Nous notons $\lambda=\lambda^\emptyset_\emptyset$ l'invariant $\lambda$ correspondant.\smallskip

\item[(ii)] Le cas $S=\emptyset$, $T=\{\ell\}$ est celui considéré par Kuz'min: le groupe $\,\C^{\{\ell\}}(K_\%)=\C^{\{\ell\}}_{\,\emptyset}(K_\%)$ s'identifie au groupe de Galois de la pro-$\ell$-extension non-ramifiée $\ell$-décomposée maximale de $K_\%$; et le $\ell$-groupe de classes associé  $\,\Cl^{\{\ell\}}(K_\%)=\Cl^{\{\ell\}}_{\,\emptyset}(K_\%)$ n'est autre que le $\ell$-groupe des $\ell$-classes de diviseurs du corps $K_\%$. Nous notons $\lambda^{\si{\{\ell\}}}=\lambda^{\si{\{\ell\}}}_{\,\emptyset}$ l'invariant $\lambda$ correspondant.\smallskip

\item[(iii)] Le cas $S=\{\ell\}$, $T=\emptyset$ est celui considéré par Wingberg: le groupe $\,\C_{\{\ell\}}(K_\%)=\C_{\{\ell\}}^{\,\emptyset}(K_\%)$ s'identifie au groupe de Galois de la pro-$\ell$-extension $\ell$-ramifiée maximale de $K_\%$; et le $\ell$-groupe de classes associé  $\,\Cl_{\{\ell\}}(K_\%)=\Cl_{\{\ell\}}^{\,\emptyset}(K_\%)$ n'est autre que le $\ell$-groupe des $\ell$-classes de diviseurs du corps $K_\%$. Nous notons $\lambda_{\si{\{\ell\}}}=\lambda_{\si{\{\ell\}}}^{\,\emptyset}$ l'invariant $\lambda$ correspondant.\smallskip
\end{itemize}

La correspondance entre ces divers paramètres résulte d'arguments classiques de dualité qui font intervenir les groupes de racines d'ordre $\ell$-primaire de l'unité. Le corps $K$ étant supposé donné, nous sommes donc amenés à raisonner dans l'extension cyclotomique $K'=K[\mmu_\ell]$ engendrée sur $K$ par les racines $\ell$-ièmes de l'unité. Les paramètres à considérer sont ainsi ceux du corps $K'$.\smallskip

Du point de vue galoisien, la situation se présente donc comme suit:
\begin{itemize}
\item $K_\%$ est un corps surcirculaire à conjugaison complexe;
\item $K_\%^+$ est son sous-corps réel;
\item $K'$ est l'extension cyclotomique $K[\mmu_\ell]$.
\end{itemize}
Le groupe de Galois $\Delta=\Gal(K'_\%/K^+_\%)$ est donc un groupe abélien d'ordre $d$ étranger à $\ell$ et d'exposant divisant $(\ell-1)$. Cette dernière condition entraîne que l'algèbre $\ell$-adique $\Zl[\Delta]$ est une algèbre semi-locale, produit direct de $d$ exemplaires de l'anneau $\Zl$, indexés par les caractères $\ell$-adiques irréductibles du groupe $\Delta$. En particulier, les invariants respectivement attachés à $K_\%$ ou à $K_\%^+$ s'interprètent comme les $\psi$-parties de ceux attachés à $K'_\%$, pour des caractères $\ell$-adiques convenables $\psi$ de $\Delta$. Par exemple, la composante réelle de $K_\%$ correspond au caractère unité $1$; et la composante imaginaire à l'unique caractère $\bar 1$ de $\Delta$ qui relève le caractère d'augmentation de $\Gal(K_\%/K_\%^+)$.\smallskip

Dans ce contexte, il est commode de regarder les invariants lambda non plus comme des entiers mais comme des caractères du groupe $\Delta$:

\begin{Conv}\label{Conv}
Soient,comme plus haut $\ell$ un nombre premier impair, $S$ et $T$ deux ensembles finis disjoints de places de $\QQ$, puis $K_\%$ un corps surcirculaire à conjugaison complexe, $K^+_\%$ son sous-corps réel maximal, $K'_\%$ le corps $K_\%[\zeta_\ell]$; et $\Delta$ le groupe de Galois $\Gal(K'_\%/K^+_\%)$.\par
Par $\lambda^S_T$ nous entendons désormais le caractère $\ell$-adique du groupe $\Delta$ construit sur les $\varphi$-composantes du groupe $\,\C^S_T(K'_\%)$; c'est-à-dire que, pour chaque caractère $\ell$-adique irréductible $\varphi$ de $\Delta$, l'invariant numérique lambda habituel de la $\varphi$-composante du $\Zl[\Delta]$-module $\,\C^S_T(K'_\%)$ est donné par le produit scalaire $\langle\lambda^S_T,\varphi\rangle$.
\end{Conv}

Cela posé, convenons de dire qu'un caractère $\ell$-adique irréductible de $\Delta$ est réel ou imaginaire suivant qu'il agit trivialement ou non sur ma conjugaison complee $\bar\tau$; et écrivons:\smallskip

\centerline{$\chi=\chi^{\si{\oplus}}+\chi^{\si{\ominus}}$}\smallskip

\noindent la décomposition d'un caractère $\ell$-adique en ses composantes réelle et imaginaire. L'application\smallskip

\centerline{$\varphi \mapsto=\varphi^*=\omega\varphi^{-\si{1}}$,}\smallskip

\noindent où $\omega$ est le caractère cyclotomique (i.e. le caractère de l'action de $\Delta$ sur les racines de l'unité) et $\varphi^{-\si{1}}$ le contragrédient de $\varphi$ (défini par $\varphi^{-\si{1}}(\sigma)=\varphi(\sigma^{-\si{1}})$) est une involution du groupe des caractères $\ell$-adiques virtuels de $\Delta$ appelée, dans la tradition bisontine, {\em involution du miroir}.\smallskip

Avec ces notations, les résultats de dualité de \cite{J18} s'énoncent comme suit:

\begin{Prop}\label{P4}
Dés lors que le corps $K'_\%$ satisfait les conjectures de Leopoldt et de Gross-Kuz'min, regardés comme caractères du groupe $\Delta$ les trois invariants $\lambda$, $\lambda_{\{\ell\}}$ et $\lambda^{\{\ell\}}$ (attachés respectivement au $\ell$-groupe des classes au sens ordinaire $\,\Cl(K'_\%)$, au $\ell$-groupe des $\ell$-classes $\,\Cl_{\{\ell\}}(K'_\%)$ et au $\ell$-groupe des classes infinitésimales $\,\Cl^{\{\ell\}}(K'_\%)$) sont liés par les identités:
\smallskip

\centerline{$\lambda_{\{\ell\}}=\lambda-\chi_\ell^{\si{\ominus}}$\qquad et \qquad $\lambda^{\{\ell\}}=\lambda^*+(\chi_\ell^{\si{\oplus}}-1)^*$,}\smallskip

\noindent où $\chi_{\ell}^\ph=\sum_{\l^+_\%\mid\ell}\chi_{\l^+_\%}$ désigne la somme des induits à $\Delta$ des caractères unités des sous-groupes de décomposition dans $K'_\%/K^+_\%$ des places de $K^+_\%$ au-dessus de $\ell$.
\end{Prop}

Plus précisément, posant de même $\chi^\ph_S=\sum_{\p^+_\%\mid S}\chi^\ph_{\p^+_\%}$, nous avons ici:

\begin{Th}\label{T1}
Soient $S$ et $T$ deux ensembles finis disjoints de places de $\QQ$ ne contenant pas $\ell$, puis $\bar S=S\cup\{\ell\}$ et $\bar T=T\cup\{\ell\}$.
Avec les conventions précédentes et sous les conjectures de Leopodt et de Gross-Kuz'min, les paramètres lambda attachés aux $\ell$-groupes de $T$-classes $S$-infinitésimales du corps $K'_\%$, regardés comme caractères du groupe $\Delta$, vérifient les identités:\medskip
\begin{itemize}
\item[(i)] \;\;$\lambda_T^S{}^{\si{\ominus}} = \lambda^{\si{\ominus}}+\omega\,[(\chi_S^{\si{\oplus}}-1)\vee 0]$.\medskip

\item[(ii)] \;$\lambda_{\bar T}^S{}^{\si{\ominus}} = \lambda_{\{\ell\}}{}^{\si{\ominus}}+\omega\,[(\chi_S^{\si{\oplus}}-1)\vee 0] = 
(\lambda-\chi_{\ell}^\ph)^{\si{\ominus}}+\omega\,[(\chi_S^{\si{\oplus}}-1)\vee 0]$.\medskip

\item[(iii)] $\lambda_T^{\bar S}{}^{\si{\oplus}} = \lambda^{\{\ell\}}{}^{\si{\oplus}}+\omega\,\chi_S^{\si{\ominus}}=
\lambda^{*\,\si{\oplus}}+\omega\,\chi_S^{\si{\ominus}}$.
\end{itemize}
\end{Th}

\Preuve Nous avons supposé ici que ni $S$ ni $T$ ne contiennent $\ell$. Il en résulte immédiatement que les places de $T$ sont sans influence sur les paramètres lambda, puisque toute $\ell$-extension de corps surcirculaires qui est non-ramifiée au-dessus de $T$ est nécessairement $T$-décomposée.
Reste à déterminer l'influence des places de $S$, c'est-à-dire la nature, dans chacun des trois groupes de galois $\,\C^T_S(K'_\%)$, $\,\C^{\bar T}_S(K'_\%)$ et $\,\C^T_{\bar S}(K'_\%)$, du sous-module engendré par les places d'inertie attachés aux places de $S$.

Prenons appui pour cela sur la théorie du corps de classes: pour chaque place $\p'_\%$ de $K'_\%$ au-dessus de $S$ et tout $n$ assez grand, le $\ell$-sous-groupe des unités principales du complété en $\p'_\%$ du sous-corps $K'_n$ de $K_\%$ de degré $\ell^n$ sur $K'$ s'identifie au $\ell$-sous-groupe de Sylow $\mu'_n$ du groupe des racines de l'unité de $K'_n$. La limite projective $\,\U_{\p'_\%}=\varprojlim \mu'_n$ s'identifie donc au module de Tate\smallskip

\centerline{$\Tl=\varprojlim \mmu_{\ell^n}$.}\smallskip

Maintenant, le groupe de Galois $\Delta$ permute transitivement les places $\p'_\%$ de $K'_\%$ au-dessus d'une même place $\p_\%^{\si{+}}$ de $K_\%^+$; et il opère sur $\Tl$ par le caractère cyclotomique $\omega$. Le produit\smallskip

\centerline{$\U_{\p_\%^{\si{+}}}=\prod_{\p'_\%\mid\p_\%^{\si{+}}}\U_{\p'_\%}$}\smallskip

\noindent est donc un $\Zl[\Delta]$-module projectif qui a pour caractère le translaté par $\omega$ de l'induit à $\Delta$ du sous-groupe de décomposition de $\p_\%^{\si{+}}$; c'est à dire le reflet $\chi^*_{\p_\%^{\si{+}}}=\omega\chi^{-\si{1}}_{\p_\%^{\si{+}}}=\omega\chi^\ph_{\p_\%^{\si{+}}}$ du caractère $\chi^\ph_{\p_\%^{\si{+}}}$.\smallskip

$(i)$ Le sous-module de $\,\C^T_S(K'_\%)^{\si{\ominus}}$ engendré par l'image des sous-groupes d'inertie attachés aux places au-dessus de $S$ est le quotient de la composante imaginaire du produit $\,\U_S=\prod_{\p_\%^{\si{+}}\mid S}\,\U_{\p_\%^{\si{+}}}$ par l'image des pro-$\ell$-groupes d'unités des $K'_n$. Comme, restreints à leur composante imaginaire, ces groupes se réduisent aux sous-groupes $\mu'_n$, l'image de leur limite projective s'identifie au  module de Tate $\Tl=\varprojlim \mmu_{\ell^n}$. Le caractère cherché est donc donné par la différence $(\omega\chi_S^\ph)^{\si{\ominus}}-\omega$, sauf si le terme de gauche est nul, ce qui ne se produit que lorsque $S$ est vide, auquel cas le caractère cherché est nul. D'où le résultat attendu (où $\alpha\vee\beta$ désigne le plus petit caractère contenant à la fois $\alpha$ et $\beta$):  $((\omega\chi_S^\ph)^{\si{\ominus}}-\omega)\vee 0=\omega((\chi_S^{\si{\oplus}}-1)\vee 0$.\smallskip

$(ii)$ Le cas de $\,\C^{\bar T}_S(K'_\%)^{\si{\ominus}}$ se traite de la même façon, à ceci près qu'il convient de remplacer les unités par les $\ell$-unités. Mais, comme sous la conjecture de Gross-Kuz'min, la composante imaginaire de la limite projective des $\ell$-unités est identique à celle des unités (cf. \cite{J18}, Ch. IV, \S2), le caractère cherché est le même que précédemment.\smallskip

$(iii)$ Enfin, pour $\,\C_{\bar S}^TK'_\%)^{\si{\oplus}}$, ce sont les unités infinitésimales qui remplacent les unités. Sous la conjecture de Leopoldt, elles sont triviales et le caractère cherché est $(\omega\chi_S^\ph)^{\si{\oplus}}=\omega\chi_S^{\si{\ominus}}$.

\section{La formule des genres dans une $\ell$-extension de corps surcirculaires}

\subsection{Position du problème}

Nous avons vu dans la section 1 que le choix de deux ensembles finis $S$ et $T$ de places de $\QQ$ permet d'attacher, sous certaines hypothèses,  à chaque corps surcirculaire à conjugaison complexe $L_\%$ un $\Zl$-module libre, disons $\X_{L_\%}$, de dimension finie, et un $\Zl$-module divisible de même codimension\footnote{la codimension d'un $\Zl$-module divisible est la dimension sur $\Zl$ de son dual de Pontrjagin.} $\X'_{L_\%}$; à savoir la composante réelle ou imaginaire du groupe de Galois $\,\C_S^T(L_\%)$ de la pro-$\ell$-extension abélienne $T$-décomposée $S$-ramifiée maximale de $L_\%$;  et celle du sous-module de $\Zl$-torsion du $\ell$-groupe $\,\Cl_S^T(L_\%)$  des $T$-classes $S$-infinitésimales de  $L_\%$.\smallskip

Ces deux groupes sont canoniquement des modules galoisiens, de sorte que si $L_\%$ est une $\ell$-extension d'un sous-corps surcirculaire à conjugaison complexe $K_\%$, le groupe $G=\Gal(L_\%/K_\%)$ opère de façon naturelle sur $\X_{L_\%}$ comme sur $\X'_{L_\%}$. Nous nous proposons de déterminer le caractère de cette action. Nous allons voir qu'il suffit pour ce faire de savoir le calculer dans un cas très particulier.\smallskip

L'argument essentiel repose sur l'existence de deux morphismes naturels entre les groupes $\X_{L_\%}$ et $\X_{K_\%}$ d'une part, $\X'_{L_\%}$ et $\X'_{K_\%}$ d'autre part: le premier $j_{\si{L_\%/K_\%}}$ (respectivement $j'_{\si{L_\%/K_\%}}$), induit par l'extension des diviseurs, correspond au transfert en termes galoisiens; le second $N_{\si{L_\%/K_\%}}$ (respectivement $N'_{\si{L_\%/K_\%}}$), qui provient de la norme arithmétique, correspond à la restriction en termes de groupes de Galois. De plus, les identités entre opérateurs\smallskip

\centerline{$N_{\si{L_\%/K_\%}}\si{\circ}\, j_{\si{L_\%/K_\%}} =[L_\%:K_\%]$ \qquad et \qquad $N'_{\si{L_\%/K_\%}}\si{\circ}\, j'_{\si{L_\%/K_\%}} =[L_\%:K_\%]$}\smallskip

\noindent montrent que $j_{\si{L_\%/K_\%}}$ (resp. $j'_{\si{L_\%/K_\%}}$ est un pseudo-isomorphisme de $\X_{K_\%}$ vers $\X_{L_\%}^G$ (resp. de $\X'_{K_\%}$ vers $\X'{}_{\!\!L_\%}^G$). Les identités qui s'ensuivent\smallskip  

\centerline{$\dim_\Zl\X^G_{L_\%}=\dim_\Zl \X^\ph_{K_\%}$ \qquad et \qquad$\codim_\Zl \X'{}_{\!\!L_\%}^G=\codim_\Zl \X'^\ph_{K_\%}$}\smallskip

\noindent vont nous permettre de dévisser la structure de $\X_{L_\%}$ comme de $\X'_{L_\%}$. Distinguons trois étapes:\medskip

\noindent{\bf Étape 1:} réduction au cas cyclique\smallskip

Pour établir que deux représentations galoisiennes sont isomorphes, il suffit de vérifier qu'elles ont même caractère. Or, l'égalité de deux caractères se teste sur les éléments du groupe de Galois. Et la valeur d'un caractère sur un élément est celle de la restriction de ce caractère au sous-groupe cyclique engendré par cet élément. Nous ne restreignons donc pas la généralité de la démonstration, si nous supposons que $G$ est un $\ell$-groupe cyclique.\medskip

\noindent{\bf Étape 2:} réduction au calcul du rang\smallskip

Si $G$ est cyclique d'ordre $\ell^m$, la décomposition de l'algèbre de groupe $\Ql[G]$ s'écrit:\smallskip

\centerline{$\Ql[G]\simeq\bigoplus_{i=0}^m\Ql[\zeta_{\ell^i}]$;}\smallskip

\noindent et celle du caractère régulier est donnée par:\smallskip

\centerline{$\Reg_G=\sum_{i=0}^m\psi_i$, avec $\psi_0=1_G$ et $\deg\psi_i=(\ell-1)\ell^{i-\si{1}}$ pour $i=1,\dots,m$.}\smallskip

\noindent En particulier, tout caractère $\ell$-adique $\chi$ du groupe $G$ est alors de la forme \smallskip

\centerline{$\chi=\sum_{i=0}^m n_i\psi_i$, avec $(n_{\si{0}},\dots,n_m)\in\NN^{m+{\si1}}$;}\smallskip

\noindent c'est-à-dire qu'il est entièrement déterminé par la donnée des $m+1$ entiers $n_{\si{0}},\dots,n_m$.\smallskip

Cela étant, si $M$ est un $\Ql[G]$-module noethérien, et si $M_0,\dots,M_m$ sont les sous-modules fixés par les $m+1$ sous-groupes $G_i$ de $G$ d'ordres respectifs $\ell^i$, les relations\smallskip

\centerline{$\dim_\Ql M_k=\sum_{i=1}^k n_i \deg \psi_i$}\smallskip

\noindent montrent que la détermination du caractère de $M$ et le calcul des dimensions des $M_i$ sont un seul et même problème.\medskip

\noindent{\bf Étape 3:} réduction au cas cyclique élémentaire\smallskip

Ce point acquis, il est clair maintenant qu'il nous suffit de donner une démonstration dans le seul cas où le groupe $G$ est cyclique d'ordre premier. En effet, tout $\ell$-groupe étant résoluble, une $\ell$-extension (galoisienne) de corps de nombres s'obtient toujours par empilement d'un nombre fini de $\ell$-extensions cycliques élémentaires. 
Si donc nous venons à disposer d'une formule de rang dans le cas cyclique de degré premier, nous tiendrons de ce fait une formule de rang valable dans toute $\ell$-extension; et donc, en fin de compte, une identité sur les caractères dans une telle extension.\smallskip

\noindent{\bf Conclusion:} schéma de démonstration\smallskip

Dans tout ce qui suit, nous supposons que $G$ est un $\ell$-groupe cyclique de degré premier (impair) et nous raisonnons directement sur les caractères des modules considérés. D'après ce qui précède, les formules obtenues dans ce cas particulier valent en fait pour un $\ell$-groupe quelconque. Les équations aux dimensions s'en déduisent simplement en prenant les degrés des caractères.
 
\subsection{Étude algébrique du cas cyclique de degré premier}

Supposons donc le groupe $G$ cyclique d'ordre $\ell$; notons $\Reg_G$ le caractère de la représentation régulière et $\Aug_G=\Reg_G-1_G$ celui de la représentation d'augmentation.

Distinguons deux cas, suivant que le $\Zl[G]$-module étudié $\X$ est libre ou divisible.


\subsubsection{Modules libres}

Si $\X$ est $\Zl$-libre et de type fini, sa décomposition comme somme directe de $\Zl[G]$-modules indécomposables est de la forme:

\centerline{$\X\simeq\Zl[G]^\alpha \oplus \Zl[\zeta_\ell]^\beta \oplus \Zl^\gamma$,}\smallskip

\noindent avec $\alpha$, $\beta$ et $\gamma$ dans $\NN$. Le caractère qui lui correspond est donc:\smallskip

\centerline{$\chi=\alpha\Reg_G + \beta \Aug_G + \gamma 1_G =(\alpha+\gamma)\Reg_G-(\gamma-\beta)\Aug_G$;}\smallskip

\noindent et la dimension $x$ de $\X$ est ainsi: $x=\ell(\alpha+\gamma)-(\ell-1)((\gamma-\beta)$.\smallskip

Dans l'identité obtenue, les deux indices $\alpha+\gamma$ et $\gamma-\beta$ ont une interprétation simple:\smallskip
\begin{itemize}
\item Le premier est la dimension du sous-module des points fixes: $\alpha+\gamma=\dim_\Zl \X^G$.\smallskip

\item Le second n'est autre que le quotient de Herbrand dimensionnel $q(\X)$ du module $\X$:\smallskip

\centerline{$\gamma-\beta =h^2(\X)-h^1(\X)$ avec $\gamma=h^2(\X)=\dim_\Fl H^2(G,\X)$ et $\beta=h^1(\X)=\dim_\Fl H^1(G,\X)$.}
\end{itemize}\smallskip

Dans ce cas, nous obtenons donc:

\begin{Lem}\label{L2}
Si $G$ est un groupe cyclique d'ordre $\ell$, le caractère d'un $\Zl[G]$-module $\X$ noethérien et sans $\Zl$-torsion est donné par la formule:

\centerline{$\chi=\dim_\Zl \X^G \Reg_G -\,q(\X)\,\Aug_G$,}\smallskip

\noindent où $q(\X)$ représente le quotient de Herbrand dimensionnel de $\X$. Et $x=\dim_\Zl X$ est ainsi:\smallskip

\centerline{$x=\ell\dim_\Zl \X^G - (\ell-1)q(\X)$.}
\end{Lem}

\subsubsection{Modules divisibles}

Si $\X'$ est $\Zl$-divisible et de cotype fini, sa décomposition comme somme directe de $\Zl[G]$-modules indécomposables est de la forme:

\centerline{$X'\simeq (\Ql/\Zl)\otimes_\Zl\Zl[G]^{\alpha'} \oplus (\Ql/\Zl)\otimes_\Zl\Zl[\zeta_\ell]^{\beta'} \oplus (\Ql/\Zl)^{\gamma'}$,}\smallskip

\noindent avec $\alpha'$, $\beta'$ et $\gamma'$ dans $\NN$. Le caractère qui lui correspond est donc:\smallskip

\centerline{$\chi'=\alpha'\Reg_G + \beta' \Aug_G + \gamma' 1_G =(\alpha'+\gamma')\Reg_G+(\beta'-\gamma')\Aug_G$;}\smallskip

\noindent et la codimension $x'$ de $\X'$ est ainsi: $x=\ell(\alpha'+\gamma')+(\ell-1)(\beta'-\gamma')$.\smallskip

Comme plus haut, les deux indices $\alpha'+\gamma'$ et $\beta'-\gamma'$ ont une interprétation simple:\smallskip
\begin{itemize}
\item Le premier est la codimension du sous-module des points fixes: $\alpha'+\gamma'=\codim_\Zl \X'{}^G$.\smallskip

\item Le second n'est autre que le quotient de Herbrand dimensionnel $q(\X')$ du module $\X'$:\smallskip

\centerline{$\beta'-\gamma' =h^2(\X')-h^1(\X')$ avec $\beta'=h^2(\X')=\dim_\Fl H^2(G,\X')$ et $\gamma'=h^1(\X')=\dim_\Fl H^1(G,\X')$.}
\end{itemize}\smallskip

Dans ce cas, nous obtenons donc:

\begin{Lem}\label{L3}
Si $G$ est un groupe cyclique d'ordre $\ell$, le caractère d'un $\Zl[G]$-module $\Zl$-divisible et de cotype fini $\X'$ est donné par la formule:

\centerline{$\chi'=\codim_\Zl \X'{}^G \Reg_G +\,q(\X')\,\Aug_G$,}\smallskip

\noindent où $q(\X')$ représente le quotient de Herbrand dimensionnel de $\X'$. Et $x'=\codim_\Zl \X$ est ainsi:\smallskip

\centerline{$x'=\ell\codim_\Zl \X'{}^G + (\ell-1)q(\X')$.}
\end{Lem}

\Nota
Dans les applications arithmétiques plus loin, on a toujours $x'=x$ et $q(\X')=-q(\X)$, de sorte qu'il est possible de raisonner indifféremment sur $\X$ ou sur $\X'$. Cela ne veut pas dire, en revanche, que l'on a systématiquement $\beta'=\beta$ et $\gamma'=\gamma$, le groupe $\X$ n'étant pas, en général, isomorphe au dual de Pontrjagin de $\X'$ en tant que module galoisien (cf. Scolie \ref{S} infra); l'égalité des caractères signifie seulement que cet isomorphisme a lieu après extension des scalaires à $\Ql$.

\subsubsection{Énoncé des résultats}

Dans la section précédente, nous avons convenu de regarder les paramètres $\lambda^T_S$ associés aux groupes $\,\C^T_S$ ou $\,\Cl^T_S$ attachés à un corps surcirculaire à conjugaison complexe non plus comme des entiers, mais comme des caractères $\ell$-adiques d'un groupe abélien fini $\Delta$  d'exposant divisant $\ell-1$.\smallskip

Pour reformuler dans ce contexte les Lemmes \ref{L2} et \ref{L3}, nous considérons la situation suivante:\smallskip
\begin{itemize}
\item $K_\%$ est un corps surcirculaire à conjugaison complexe; $K^+_\%$ est son sous-corps réel maximal; $K'_\%=K_\%[\mmu_\ell]$ est engendré par les racines $\ell$-ièmes de l'unité pour un premier impair $\ell$.
\item $L_\%$ est une $\ell$-extension (galoisienne) de $K_\%$ qui provient par composition avec $K_\%$ d'une $\ell$-extension totalement réelle $L^+_\%$ de $K^+_\%$; et $L'_\%$ est le corps $L_\%[\mmu_\ell]$.
\item $G$ est le groupe de Galois $\Gal(L'_\%/K'_\%)\simeq\Gal(L^+_\%/K^+_\%)$; et $\Delta$ est le groupe abélien $\Gal(L'_\%/L^+_\%)\simeq\Gal(K'_\%/K^+_\%)$.
\item Les objets étudiés sont équipés d'une action naturelle de $G$ et de $\Delta$, c'est-à-dire que ce sont des $\Zl[\Delta\times G]$-modules.
\end{itemize}\smallskip

Cela posé, par application des arguments précédents à chacune des composantes $\Delta$-isotypiques de $\X$ et de $\X'$, nous obtenons immédiatement:

\begin{Prop}\label{P5}
Soit $G$ cyclique d'ordre $\ell$, $\Reg_G$ le caractère de sa représentation régulière et $\Aug_G=\Reg_G-1_G$ celui de sa représentation d'augmentation. Alors:
\begin{itemize}
\item[(i)] Si $X$ est un $\Zl[\Delta\times G]$-module, de type fini et libre sur $\Zl$, son caractère $\chi$ est donné par la formule:

\centerline{$\chi=\chi_\Delta(\X^G)\Reg_G-q_\Delta(G,\X)\Aug_G $,}\smallskip

où $\chi_\Delta(\X^G)$ est le caractère de $\X^G$ regardé comme $\Zl[\Delta]$-module et $q_\Delta(G,\X)$ le $\Delta$-quotient de Herbrand du $G$-module $\X$, i.e. l'unique caractère virtuel du groupe $\Delta$ à valeurs dans $\Zl$ qui relève la différence $h^2_\Delta(G,\X)-h^1_\Delta(G,\X)$ des caractères respectifs des $\Fl[\Delta]$-modules  $H^2(G,\X)$ et $H^1(G,\X)$.
\smallskip

\item[(ii)] Si $X'$ est un $\Zl[\Delta\times G]$-module, de cotype fini et divisible sur $\Zl$, son caractère $\chi'$ est donné par la formule:\smallskip

\centerline{$\chi'=\chi'_\Delta(\X^G)\Reg_G+q_\Delta(G,\X')\Aug_G $,}\smallskip

où $\chi'_\Delta(\X^G)$ et $q_\Delta(G,\X')$ sont définis comme plus haut.
\end{itemize}
\end{Prop}

\subsection{Étude arithmétique du cas cyclique de degré premier impair}

Plaçons-nous dans  le contexte arithmétique de la Proposition \ref{P5} ci-dessus. Sous réserve de trivialité de l'invariant mu, le $\Zl$-module libre qui nous intéresse est la composante imaginaire $\X_{L'_\%}=\C_{L'_\%}^{\si{\ominus}}$ du groupe de Galois de la pro-$\ell$-extension abélienne non-ramifiée maximale de $L'_\%$; et le $\Zl$-module divisible qui lui correspond est la composante imaginaire  $\X'_{L'_\%}=\Cl_{L'_\%}^{\si{\ominus}}$ du $\ell$-groupe des classes de diviseurs (au sens ordinaire) du corps surcirculaire $L'_\%$.\smallskip

Comme plus haut, nous notons $\chi^\ph_{L'_\%}$ le caractère de $\X_{L'_\%}$ et $\chi'_{L'_\%}$ celui de $\X'_{L'_\%}$, l'un et l'autre regardés comme $\Zl[\Delta \times G]$-modules. Nous réservons la notation $\lambda_{L'_\%}$ pour désigner indifféremment le caractère de $\X_{L'_\%}$, regardé comme $\Zl[\Delta]$-module, ou celui de $\X'_{L'_\%}$. L'invariant lambda traditionnel d'Iwasawa n'est autre que le degré de $\lambda_{L'_\%}$; c'est aussi le degré commun de $\chi^\ph_{L'_\%}$ et de $\chi'_{L'_\%}$.

Cela étant, la Proposition \ref{P5} nous dit que pour exprimer $\chi^\ph_{L'_\%}$ et  $\chi'_{L'_\%}$ à partir de $\chi^\ph_{K'_\%} = \chi'_{K'_\%}$, il suffit de déterminer l'action de $\Delta$ sur les groupes de cohomologie $H^1(G,\C_{L'_\%}^{\si{\ominus}})$ et $H^2(G,\C_{L'_\%}^{\si{\ominus}})$, d'une part; $H^1(G,\Cl_{L'_\%}^{\si{\ominus}})$ et $H^2(G,\Cl_{L'_\%}^{\si{\ominus}})$, d'autre part.

\subsubsection{Cohomologie des groupes de classes}

\begin{Lem}\label{L4}
Soient $L_\%/K_\%$ une $\ell$-extension cyclique élémentaire de corps surcirculaires et $G$ son groupe de Galois. Alors:
\begin{itemize}
\item[(i)] Le groupe multiplicatif $L^\times_\%$ est cohomologiquement trivial; autrement dit, on a:\smallskip

\centerline{$H^1(G,L^\times_\%)=1$ \qquad et \qquad $H^2(G,L^\times_\%)=1$.}\smallskip

\item[(ii)] La cohomologie du groupe des diviseurs est donnée par les formules:\smallskip

\centerline{$H^1(G,D_{L_\%})=1$ \qquad et \qquad $H^2(G,P_{L_\%}) \simeq \Ram'_{L_\%/K_\%} \simeq \bigoplus_{\p_\%\nmid\ell}\ZZ/e_{\p_\%}\ZZ$,}\smallskip

\noindent où $\Ram'_{L_\%/K_\%}$ mesure la ramification modérée, $\p_\%$ parcourt les places étrangères à $\ell$ qui se ramifient dans $L_\%/K_\%$ et $e_{\p_\%}=\ell$ est l'indice de ramification.\smallskip

\item[(iii)] La cohomologie du groupe des diviseurs principaux est duale de celle du groupe des unités:\smallskip

\centerline{$H^1(G,P_{L_\%}) \simeq H^2(G,E_{L_\%})$ \qquad et \qquad $H^2(G,P_{L_\%}) \simeq H^1(G,E_{L_\%})$.}
\end{itemize}
\end{Lem}

\Preuve
La trivialité des groupes de cohomologie $H^1(G,L^\times)$ et $H^1(G,D_{L_\%})$ est bien connue, du moins pour les corps de nombres: c'est ce qu'il est convenu d'appeler le théorème 90 de Hilbert\footnote{Bien entendu, dans son traité de 1913, Hilbert n'énonçait pas son résultat sous forme cohomologique, mais seulement dans le cas cyclique: D. Hilbert, {\em Théorie des corps de nombres algébriques}, Hermann, Paris (1913).}. Son extension aux corps surcirculaires ne pose pas de difficulté: il suffit d'écrire $L_\%/K_\%$ comme limite inductive de sous-extensions de corps de nombres relativement cycliques et de leur appliquer le résultat de Hilbert. Le cas des groupes $H^2(G,L^\times)$ et $H^2(G,D_{L_\%})$, en revanche, est plus surprenant, puisque, pour les corps de nombres, aucun de ces deux groupes n'est fini.\smallskip

Cela dit, examinons successivement les diverses assertions du Lemme.\smallskip

$(i)$ L'identité $H^2(G,L^\times)=1$ est signalée par Iwasawa (\cite{Iw2}, Lem. 5) et vaut de fait pour tout $\ell$-groupe $G$. Lorsque $G$ est cyclique, elle dit simplement que les éléments de $K_\%^\times$ sont normes d'éléments de $L_\%^\times$ (cf. \cite{DM}, Lem. 3): $K_\%^\times=N_{\si{L_\%/K_\%}}(L_\%^\times)$.

Et pour l'établir dans ce cas, d'après le principe de Hasse, il suffit de vérifier que tout élément $x$ de $K_\%^\times$ est norme locale partout dans l'extension $L_\%/K_\%$. Les palces à l'infini ne faisant pas problème, puisque $\ell$ est pris impair, donnons une place finie $\p_\%$ de $K_\%$. Écrivons $L_\%=\bigcup_{n\in\NN}L_N$ et $K_\%=\bigcup_{n\in\NN}K_n$ comme réunions croissantes de corps de nombres avec $L_n/K_n$ cyclique de groupe $G$ et $[L_{n+\si{1}}:L_n]=[K_{n+\si{1}}:K_n]=\ell$; puis fixons $n$ assez grand pour avoir $x\in K_n$ et $\p_\%$ non-décomposée dans $K_\%/K_n$. Prenant alors les symboles de Hasse attachés aux places en dessous de $\p_\%$, nous obtenons:

\centerline{$\big(\frac{x,L_{n+\si{1}}/K_{n+\si{1}}}{\p^\ph_{n+\si{1}}}\big)=\big(\frac{N_{n+\si{1}/n}(x),L_n/K_n}{\p^\ph_n}\big)=\big(\frac{x^\ell,L_n/K_n}{\p^\ph_n}\big)=\big(\frac{x,L_n/K_n}{\p^\ph_n}\big)^\ell=1$,}\smallskip

\noindent et $x$ est norme locale en $\p^\ph_{n+\si{1}}$ dans $L_{n+\si{1}}/K_{n+\si{1}}$, donc norme locale en $\p_\%$ dans $L_\%/K_\%$.\smallskip

$(ii)$ L'identité $H^2(G,D_{L_\%}) \simeq \Ram'_{L_\%/K_\%}$, (où le groupe $\Ram'_{L_\%/K_\%}$ est défini par l'énoncé) s`obtient comme suit: $G$ étant cyclique, nous avons l'isomorphisme:\smallskip

\centerline{$H^2(G,D_{L_\%}) \simeq D_{L_\%}^G/\nu_{L_\%/K_\%}(D_{L_\%})$,}\smallskip

\noindent où $\nu_{L_\%/K_\%}$ désigne la norme algébrique $\sum_{\sigma\in G}\sigma$; ce que nous pouvons encore écrire:\smallskip

\centerline{$H^2(G,D_{L_\%}) \simeq D_{L_\%}^G/N_{L_\%/K_\%}(D_{L_\%})$,}\smallskip

\noindent où $N_{L_\%/K_\%}$ désigne cette fois la norme arithmétique, en convenant d'identifier les diviseurs de $K_\%$ avec leurs étendus à $L_\%$.  
Maintenant, tout diviseur de $K_\%$ est banalement la norme d'un diviseur de $L_\%$: les diviseurs premiers au-dessus de $\ell$, parce qu'ils sont $\ell$-divisibles; ceux étrangers à $\ell$, parce qu'ils ne présentent plus d'inertie. Il vient donc:\smallskip

\centerline{$H^2(G,D_{L_\%}) \simeq D_{L_\%}^G/(D_{K_\%})$;}\smallskip

\noindent et la formule annoncée en résulte, les diviseurs ambiges étant engendrés par les diviseurs ramifiés modulo ceux étendus. Comme les diviseurs au-dessus de $\ell$ sont $\ell$-divisibles, donc étendus, seuls restent les diviseurs ramifiés étrangers à $\ell$, ce qui donne bien:\smallskip

\centerline{$H^2(G,D_{L_\%}) \simeq \sum_{\p_\%\nmid\ell} D_{L_\%}(\p_\%)^G/(D_{K_\%}(\p_\%) \simeq \sum_{\p_\%\nmid\ell} \ZZ/e_{\p_\%}\ZZ$.}\smallskip

$(iii)$ Enfin, le fait que les groupes d'unités $E_{L_\%}$ et de diviseurs principaux $P_{L_\%}$ soient en dualité cohomologique résulte immédiatement du point $(i)$. La suite exacte courte canonique\smallskip

\centerline{$1 \to E_{L_\%} \to L_\%^\times \to P_{L_\%} \to 1$}\smallskip

\noindent fournit, en effet, l'hexagone exact, où les sommets supérieur et inférieur sont nuls:
\begin{displaymath}
\xymatrix{
 & H^1(G,L^\times_\%) \ar@{->}[dr] &\\
 H ^1(G,E_{L_\%}) \ar@{->}[ur] &&  H ^1(G,P_{L_\%}) \ar@{->}[d]^\simeq\\
 H ^2(G,P_{L_\%}) \ar@{->}[u]^\simeq &&  H ^2(G,E_{L_\%}) \ar@{->}[dl]\\ 
 & H^2(G,L^\times_\%) \ar@{->}[ul] &
}
\end{displaymath}
D'où le résultat annoncé.

\begin{Lem}\label{L5}
Supposons maintenant que $L_\%/K_\%$ est une $\ell$-extension cyclique élémentaire de corps surcirculaires à conjugaison complexe; et notons : $\Delta=\Gal(L'_\%/L^+_\%)\simeq Gal(K'_\%/K^+_\%)$. Dans ce cadre, les identités précédentes se traduisent par les isomorphismes de $\Fl[\Delta]$-modules:\smallskip
\begin{itemize}
\item[(i)] $H^1(G,P^{\,\si{\ominus}}_{L'_\%}) \simeq H^2(G,E^{\,\si{\ominus}}_{L'_\%})=1$; et  $H^2(G,P^{\,\si{\ominus}}_{L'_\%}) \simeq H^1(G,E^{\,\si{\ominus}}_{L'_\%}) \simeq \mmu_\ell$.\smallskip

\item[(ii)]  $H^1(G,D^{\,\si{\ominus}}_{L'_\%})=1$; et  $H^2(G,D^{\,\si{\ominus}}_{L'_\%}) \simeq \bigoplus_{\p_\%^{\si{+}}\nmid\ell}(\ZZ/e_{\p^+_\%}\ZZ)[\Delta/\Delta_{\p_\%^{\si{+}}}]$.\smallskip
\end{itemize}
En particulier, les caractères des groupes  $H^1(G,P_{L'_\%})$ et $H^2(G,P_{L'_\%})$ sont donnés par:\smallskip

\centerline{$h^1(G,P_{L'_\%})=\omega$ \qquad et \qquad $h^2(G,P_{L'_\%})=\sum_{\p_\%^{\si{+}}\nmid\ell}\chi_{\p_\%^{\si{+}}}^{\si{\ominus}}$,}\smallskip

\noindent Dans ces formules  $\mmu_\ell$ est le groupe des racines $\ell$--ièmes de l'unité et $\mmu_{\ell^\%}$ le groupe des racines d'ordre $\ell$-primaire; $\omega$ désigne le caractère cyclotomique; la somme à droite porte sur les places $\p_\%^{\si{+}}\nmid\ell$ de $K^+_\%$ qui se ramifient dans $L^+_\%/K_\%^+$; et $\chi_{\p_\%^{\si{+}}}^{\si{\ominus}}$ est la composante imaginaire de l'induit à $\Delta$ du caractère de la représentation unité du sous-groupe de décomposition $\Delta_{\p_\%^{\si{+}}}$ de $\p_\%^{\si{+}}$ dans $K'_\%/K^+_\%$.
\end{Lem}

\Preuve En présence d'une conjugaison complexe, les unités imaginaires se réduisent aux seules racine de l'unité. Il vient donc:\smallskip

\centerline{$H^1(G,E^{\,\si{\ominus}}_{L'_\%})=H^1(G,\mmu_{\ell^\%})=\mmu_\ell$ \qquad et \qquad $H^2(G,E^{\,\si{\ominus}}_{L'_\%})=H^(G,\mmu_{\ell^\%})=\mmu_{\ell^\%}/\mmu_{\ell^\%}^\ell=1$.}\smallskip

\noindent D'où le résultat annoncé. Les autres assertions sont immédiates en vertu du Lemme \ref{L4}.\medskip

Ce point acquis, nous pouvons énoncer:

\begin{Prop}\label{P6}
Dans une $\ell$-extension cyclique élémentaire de corps surcirculaires à conjugaison complexe le quotient de Herbrand dimensionnel du $\ell$-groupe des classes imaginaires est donné, avec les notations précédentes, par la formule:\smallskip

\centerline{$q_\Delta(G,\Cl^{\,\si{\ominus}}_{L'_\%})=h^2_\Delta(G,D^{\,\si{\ominus}}_{L'_\%})-h^2_\Delta(G,P^{\,\si{\ominus}}_{L'_\%})=\sum_{\p_\%^{\si{+}}\nmid\ell}\chi_{\p_\%^{\si{+}}}^{\si{\ominus}} -\omega$.}
\end{Prop}

\Preuve La suite exacte courte $1 \to P^{\,\si{\ominus}}_{L'_\%} \to D^{\,\si{\ominus}}_{L'_\%} \to \Cl^{\,\si{\ominus}}_{L'_\%} \to 1$ qui définit le groupe $\,\Cl^{\,\si{\ominus}}_{L'_\%}$ donne la suite exacte longue de cohomologie, où les termes extrêmes sont nuls en vertu du Lemme \ref{L5}:\smallskip

\centerline{$H^1(G,D^{\,\si{\ominus}}_{L'_\%}) \to H^1(G,\Cl^{\,\si{\ominus}}_{L'_\%}) \to H^2(G,P^{\,\si{\ominus}}_{L'_\%}) \to H^2(G,D^{\,\si{\ominus}}_{L'_\%}) \to H^2(G,\Cl^{\,\si{\ominus}}_{L'_\%}) \to H^1(G,P^{\,\si{\ominus}}_{L'_\%})$.}\smallskip

\subsubsection{Cohomologie des groupe de Galois}

L'étude des groupes de Galois est en soi inutile au regard du résultat que nous avons en vue, puisque nous savons déjà que les caractères $q_\Delta(G,\C^{\,\si{\ominus}}_{L'_\%})$ et $q_\Delta(G,\Cl^{\,\si{\ominus}}_{L'_\%})$ sont opposés en vertu de la Proposition \ref{P5}, les deux groupes $\,\C_{L'_\%}^{\,\si{\ominus}}$ et $\,\Cl_{L'_\%}^{\,\si{\ominus}}$ d'une part, $\,\C_{K'_\%}^{\,\si{\ominus}}$ et $\,\Cl_{K'_\%}^{\,\si{\ominus}}$ d'autre part, regardés comme $\Zl[\Delta]$-modules, ayant même caractère (sous la condition $\mu=0$).\smallskip

Néanmoins, outre qu'elle fournit une autre démonstration de la Proposition \ref{P6}, la description de leur cohomologie ouvre des perspectives intéressantes sur l’arithmétique des corps surcirculaires.\smallskip

Soient donc, comme plus haut, $L_\%/K_\%$ une $\ell$ extension cyclique de degré $\ell$ de corps surcirculaires à conjugaison complexe, $G$ son groupe de Galois et $\sigma$ un générateur de $G$. Conservons les conventions précédentes; introduisons le module de Tate $\Tl=\varprojlim \mmu_{\ell^n}$, qui est un $\Zl[\Delta]$-module projectif de caractère $\omega$; et faisons choix d'un générateur topologique $\zeta$ de $\Tl$, i.e. d'une famille cohérente $(\zeta_{\ell^n})_n$ de racines $\ell^n$-ièmes primitives de l'unité.\smallskip

Regardons d'abord le groupe de cohomologie $H^1(G,\C^{\,\si{\ominus}}_{L'_\%})$. Le groupe $\,\C^{\,\si{\ominus}}_{L'_\%}$ étant sans $\Zl$-torsion (sous réserve de nullité de l'invariant mu), le noyau dans $\,\C^{\,\si{\ominus}}_{L'_\%}$ de la norme algébrique $\nu_{L'_\%/K'_\%}=\sum_{i=0}^{\ell-1}\sigma^i$ est aussi celui de la norme arithmétique $N_{L'_\%/K'_\%}$, ce qui nous donne canoniquement:\smallskip

\centerline{$H^1(G,\C^{\,\si{\ominus}}_{L'_\%})={}_\nu\C^{\,\si{\ominus}}_{L'_\%}/\C^{\,\si{\ominus}\,(\sigma-\si{1})}_{L'_\%}={}_N\C^{\,\si{\ominus}}_{L'_\%}/\C^{\,\si{\ominus}\,(\sigma-\si{1})}_{L'_\%} = \big({}_N\C^\ph_{L'_\%}/\C^{\,(\sigma-\si{1})}_{L'_\%}\big)^{\si{\ominus}} $.}\smallskip

\noindent Ce point acquis, le quotient ${}_N\C^\ph_{L'_\%}/\C^{\,(\sigma-\si{1})}_{L'_\%}$ a une interprétation arithmétique simple: désignons par $\widehat K'_\%$ (respectivement $\widehat L'_\%$ la pro-$\ell$-extension abélienne non-ramifiée maximale de $K'_\%$ (resp. de $L'_\%$); et notons $\wi L'_\%$ la sous-extension maximale de $\widehat L'_\%$ qui est abélienne sur $K'_\%$, i.e. le $\ell$-corps des genres de $L'_\%/K'_\%$. Nous avons alors le diagramme, où sont représentés les corps et les groupes de Galois:

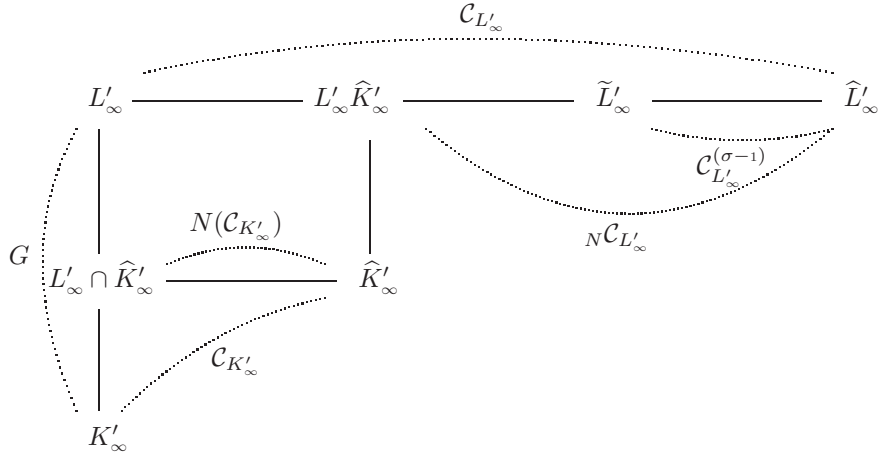
\begin{center}
\unitlength=1.5cm
\begin{picture}(8.3,4)

\put(0.7,0){$K'_\%$}
\put(0.8,0.3){\line(0,1){0.9}}
\put(0.35,1.4){$L'_\%\cap\widehat K'_\%$}
\put(0.8,1.7){\line(0,1){1.1}}
\put(0.7,3){$L'_\%$}

\bezier{50}(0.6,0.3)(0,1.7)(0.6,2.8)
\put(0,1.6){$G$}

\put(1.4,1.45){\line(1,0){1.5}}
\put(3.1,1.4){$\widehat K'_\%$}
\put(3.2,1.7){\line(0,1){1.}}
\put(1.1,3.05){\line(1,0){1.5}}
\put(2.7,3){$L'_\% \widehat K'_\%$}
\put(3.5,3.05){\line(1,0){1.5}}
\put(5.2,3){$\wi L'_\%$}

\put(5.7,3.05){\line(1,0){1.5}}
\put(7.4,3){$\widehat L'_\%$}

\bezier{120}(1.2,3.3)(4,3.9)(7.3,3.3)
\put(4.0,3.75){$\C_{L'_\%}$}

\bezier{40}(5.7,2.8)(6.5,2.6)(7.3,2.8)
\put(6.1,2.4){$\C_{L'_\%}^{(\sigma-\si{1})}$}

\bezier{100}(3.7,2.8)(5.5,1.3)(7.3,2.8)
\put(5.1,1.8){${}_N\C_{L'_\%}$}

\bezier{40}(1.4,1.6)(2.1,1.9)(2.8,1.6)
\put(1.6,1.9){$N(\C_{K'_\%})$}

\bezier{50}(1.0,0.3)(1.8,1.1)(2.8,1.3)
\put(1.8,0.7){$\C_{K'_\%}$}

\end{picture}
\end{center}\medskip

\noindent On voit par là que ${}_N\C^\ph_{L'_\%}/\C^{\,(\sigma-\si{1})}_{L'_\%}$ s'identifie au groupe de Galois de l'extension abélienne $\wi L'_\%/L'_\%\widehat K'_\%$.

Maintenant, il est bien connu, du moins pour es corps de nombres (cf. \cite{J18}, Ch.III, \S2), que ce dernier groupe est donné par la suite exacte des genres: plaçons-nous provisoirement à un étage fini $L'_n/K'_n$ de l'extension $L'_\%/K'_\%$ et définissons $\widehat K'_n$, $\widehat L'_n$ et $\wi L'_n$ comme ci-dessus. La théorie des genres nous donne alors la suite exacte courte de $\Zl[\Delta]$-modules:\smallskip

\centerline{$1 \to E_{K'_n}/(E_{K'_n}\cap N_{\si{L'_n/K'_n}}(L'_n{\!}^\times ) \overset{\h_n}{\longrightarrow} \wi\bigoplus_{\p'_n}\In_{\p'_n}(\wi L'_n/K'_n) \to \Gal(\wi L'_n/L'_n\widehat K'_n) \to 1$,}\smallskip

\noindent où $E_{K'_n}$ est le groupe des unités de $K'_n$; l'application $\h_n$ est la famille des symboles de Hasse $\big(\frac{\dot,\wi L'_n/K'_n}{\p'_n}\big)$ relatifs à l'extension abélienne $\wi L'_n/K'_n$; et le terme médian est la somme directe des sous-groupes d'inertie des places de $K'_\%$ restreinte aux familles $(\sigma_{\p'_n})_{\p'_n}$ qui vérifient la formule du produit:

\centerline{$\prod_{\p'_n}\sigma_{\p'_n}=1$.}\smallskip

Par restriction aux composantes imaginaires puis passage à la limite projective, il vient ainsi:

\begin{Lem}\label{L6}
Avec les notations ci-dessus et sous réserve de nullité de l'invariant mu, le groupe de cohomologie $H^1(G,\C^{\,\si{\ominus}}_{L'_\%})$ est donné, comme $\Zl[\Delta]$-module, par l'isomorphisme des genres:\smallskip

\centerline{$H^1(G,\C^{\,\si{\ominus}}_{L'_\%})\simeq [\oplus \In_{\p'_\%}(\wi L'_\%/K'_\%)]^{\si{\ominus}}/\h_\%(\Tl)$,}\smallskip

\noindent où la somme (restreinte à la formule du produit) porte sur les groupes d'inertie des places de $K'_\%$ dans l'extension abélienne $\wi L'_\%/K'_\%$ et l'application $\h$ est induite par les symboles de Hasse.\par
En particulier, le caractère de  $H^1(G,\C^{\,\si{\ominus}}_{L'_\%})$ est ainsi:\smallskip

\centerline{$h^1_\Delta(G,\C^{\,\si{\ominus}}_{L'_\%})=\sum_{\p_\%^{\si{+}}}\chi_{\p_\%^{\si{+}}}^{\si{\ominus}}-\delta\omega$,}\smallskip

\noindent $\chi_{\p_\%^{\si{+}}}^{\si{\ominus}}$ et $\omega$ étant définis comme dans le Lemme \ref{L5} et l'entier $\delta$ valant $0$ ou $1$ suivant que le générateur topologique $\zeta$ de $\Tl$ est norme ou pas dans $L'_\%/K'_\%$.
\end{Lem}

\Preuve Par restriction aux composantes imaginaires, la suite exacte des genres prend la forme:\smallskip

\centerline{$1 \to \mu_{K'_n}/(\mu_{K'_n}\cap N_{\si{L'_n/K'_n}}(L'_n{\!}^\times ) \overset{\h_n}{\longrightarrow} [\wi\bigoplus_{\p'_n}\In_{\p'_n}(\wi L'_n/K'_n)]^{\si{\ominus}} \to \Gal(\wi L'_n/L'_n\widehat K'_n)^{\si{\ominus}} \to 1$;}\smallskip

\noindent ce qui nous donne, par passage à la limite projective, la suite exacte de $\Zl[\Delta]$-modules finis:\smallskip

\centerline{$1 \to \Tl/\Tl^{\ell^\delta} \overset{\h_\%}{\longrightarrow} [\bigoplus_{\p'_\%}\In_{\p'_\%}(\wi L'_\%/K'_\%)]^{\si{\ominus}} \to \Gal(\wi L'_\%/L'_\%\widehat K'_\%)^{\si{\ominus}} \to 1$;}\smallskip

\noindent d'où la première assertion du lemme. Notons que la formule du produit, qui est réelle, disparait côté imaginaire. L'expression du caractère $h^1_\Delta(G,\C_{L'_\%})$ s'en déduit immédiatement, via l'isomorphisme:\smallskip

\centerline{$[\bigoplus_{\p'_\%}\In_{\p'_\%}(\wi L'_\%/K'_\%)]^{\si{\ominus}} \simeq [\bigoplus_{\p'_\%}\In_{\p'_\%}(L'_\%/K'_\%)]^{\si{\ominus}} \simeq [\bigoplus_{\p^{\si{+}}_\%}\In_{\p'_\%}(\wi L^+_\%/K^+_\%)]^{\si{\ominus}} $.}

\begin{Lem}\label{L7}
Avec les notations du Lemme \ref{L6} le caractère du groupe $H^2(G,\C_{L'_\%}^{\,\si{\ominus}})$ est donné par la formule (où la somme porte sur les places de $K^+_\%$ au-dessus de $\ell$ qui se ramifient dans $L_\%^+/K_\%^+$):\smallskip

\centerline{$h^2_\Delta(G,\C_{L'_\%}^{\,\si{\ominus}})= \sum_{\p_\%^{\si{+}}\mid\ell}\chi_{\p_\%^{\si{+}}}^{\si{\ominus}}+(1-\delta)\omega$.}
\end{Lem}

\Preuve
Regardons le groupe de cohomologie $H^2(G,\C_{L'_\%}^{\,\si{\ominus}})\simeq\C_{L'_\%}^{\,\si{\ominus}\,G}/\C_{L'_\%}^{\,\si{\ominus}\,\nu}$ comme limite projective des groupes $H^2(G,\Cl_{L'_n}^{\,\si{\ominus}})\simeq\Cl_{L'_n}^{\,\si{\ominus}\,G}/\Cl_{L'_n}^{\,\si{\ominus}\,\nu}$. Dans ce dernier quotient $\,\Cl_{L'_n}^{\,\si{\ominus}\,G}$ est la composante imaginaire du sous-groupe ambige de $\,\Cl_{L'_n}$; et $\,\Cl_{L'_n}^{\,\si{\ominus}\,\nu}=j_{L'_n/K'_n}\si{\circ}N_{\si{L'_n/K'_n}}(\Cl_{L'_n}^{\,\si{\ominus}})=\Cl_{L'_n}^{\,\et\,\si{\ominus}}$ celle du sous-groupe de $\,\Cl_{L'_n}^{\,\si{\ominus}}$ engendré par les classes des étendus à $L'_n$ des diviseurs imaginaires de $K'_n$\footnote{Le quotient $\,\Cl_{K'_n}/N_{\si{L'_n/K'_n}}(\Cl_{L'_n})\simeq\Cl_{K^+_n}/N_{\si{L^+_n/K^+_n}}(\Cl_{L^+_n}) $ étant réel, en vertu de l'isomorphisme du corps de classes, la restriction aux composantes imaginaires de la norme arithmétique est un épimorphisme.}.

Introduisons le groupe intermédiaire $\,\Cl_{L'_n}^{\,\amb\,\si{\ominus}}$ construit sur les classes des diviseurs imaginaires ambiges (i.e. des diviseurs imaginaires invariants par $G$); et formons la suite exacte courte (qui est une étape classique de la démonstration par Chevalley de sa formule des classes ambiges):\smallskip

\centerline{$1 \to \Cl_{L'_n}^{\,\amb\,\si{\ominus}}/\Cl_{L'_n}^{\,\et\,\si{\ominus}} \to \Cl_{L'_n}^{\,G\,\si{\ominus}}/\Cl_{L'_n}^{\,\et\,\si{\ominus}} \to \Cl_{L'_n}^{\,G\,\si{\ominus}}/\Cl_{L'_n}^{\,\amb\,\si{\ominus}} \to 1$.}\smallskip

Dans celle-ci le groupe à droite s'identifie au quotient\footnote{L'isomorphisme annoncé s'obtient en envoyant la classe d'un diviseur imaginaire $\a$ représentant une classe de $\,\Cl^{\,G}_{L'_n}$ sur l'unité $\varepsilon=N_{\si{L'_n/K'_n}}(\alpha)$, où $\alpha$  est un générateur imaginaire du diviseur principal $\a^{\sigma-\si{1}}$.}\smallskip

\centerline{$\,\Cl_{L'_n}^{\,G\,\si{\ominus}}/\Cl_{L'_n}^{\,\amb\,\si{\ominus}} \simeq \big(E_{K'_n}\cap N_{\si{L'_n/K'_n}}(L'_n{\!}^\times )/N_{\si{L'_n/K'_n}}(E_{L'_n})\big)^{\si{\ominus}} \simeq \mu_{K'_n}\cap N_{\si{L'_n/K'_n}}(L'_n{\!}^\times )/\mu_{K'_n}^\ell$.}\smallskip

\noindent Et, par passage à la limite projective, nous obtenons donc, avec des notations évidentes, la suite:\smallskip

\centerline{$1 \to \C_{L'_\%}^{\,\amb\,\si{\ominus}}/\C_{L'_\%}^{\,\et\,\si{\ominus}} \to \C_{L'_\%}^{\,G\,\si{\ominus}}/\C_{L'_\%}^{\,\et\,\si{\ominus}} \to \C_{L'_\%}^{\,G\,\si{\ominus}}/\C_{L'_\%}^{\,\amb\,\si{\ominus}} \to 1$,}\smallskip

\noindent où le quotient à droite fait ainsi apparaître l'entier $\delta$ défini plus haut: $\,\C_{L'_\%}^{\,G\,\si{\ominus}}/\C_{L'_\%}^{\,\amb\,\si{\ominus}} \simeq \Tl^{\ell^\delta}/\Tl^\ell$.

Reste donc à examiner le quotient à gauche. Aux niveaux finis les diviseurs imaginaires ambiges de $L'_n$ sont engendrés par les diviseurs imaginaires étendus de $K'_n$ et par ceux ramifiés dans $L'_n/K'_n$. Les ramifiés étrangers à $\ell$ sont ultimement inertes dans $L'_\%/L'_n$ donc disparaissent dans la limite projective. Le quotient $\,\C_{L'_\%}^{\,\amb\,\si{\ominus}}/\C_{L'_\%}^{\,\et\,\si{\ominus}}$ est ainsi engendré par les classes des familles cohérentes construites sur les seules places au-dessus de $\ell$ qui sont décomposées dans $L/L^+$ et ramifiées dans $L/K$.Ce point acquis, fixons $n_{\si{0}}$ assez grand pour que les diviseurs imaginaires primitifs (i.e. non puissances $\ell$-ièmes) de $L'_{n_{\si{0}}}$ construits sur les places au-dessus de $\ell$ ne soient pas principaux\footnote{L'existence d'un tel $n_{\si{0}}$ fait l'objet de la Proposition \ref{P7} infra.}.
À chaque tel diviseur $\a_{n_{\si{0}}}$ correspond alors une classe d'ordre $\ell$ dans le quotient $\,\C_{L'_{n_{\si{0}}}}^{\,\amb\,\si{\ominus}}/\C_{L'_{n_{\si{0}}}}^{\,\et\,\si{\ominus}}$ et plus généralement dans chacun des quotients $\,\C_{L'_n}^{\,\amb\,\si{\ominus}}/\C_{L'_n}^{\,\amb\,\si{\ominus}}$ pour $n\geqslant n_{\si{0}}$, à savoir la classe du diviseur $\a_n =\a _{n_{\si{0}}} ^{1/\ell^ {n-n_{\si{0}}} }$, donc finalement une classe d'ordre $\ell$ dans la limite projective $\, \C_{L'_\%}^{\,\amb\,\si{\ominus}}/\C_{L'_\%}^{\,\et\,\si{\ominus}}$.

Il suit de là que le quotient étudié $\, \C_{L'_\%}^{\,\amb\,\si{\ominus}}/\C_{L'_\%}^{\,\et\,\si{\ominus}}$ est un $\Fl[\Delta]$-module de caractère
$\sum_{\p_\%^{\si{+}}\mid\ell}\chi_{\p_\%^{\si{+}}}^{\si{\ominus}}$.\medskip

Réunissant les Lemmes \ref{L6} et \ref{L7}, nous obtenons une seconde version de la Proposition \ref{P6}:

\begin{Prop}\label{P6'}
Dans une $\ell$-extension cyclique élémentaire de corps surcirculaires à conjugaison complexe le quotient de Herbrand dimensionnel de la composante imaginaire du groupe de Galois $\,\C_{L'_\%}$ est donnée, sous réserve de nullité de l'invariant mu, par la formule:\smallskip

\centerline{$q_\Delta(G,\C^{\,\si{\ominus}}_{L'_\%})=h^2_\Delta(G,\C^{\,\si{\ominus}}_{L'_\%})-h^1_\Delta(G,\C^{\,\si{\ominus}}_{L'_\%})=\omega-\sum_{\p_\%^{\si{+}}\nmid\ell}\chi_{\p_\%^{\si{+}}}^{\si{\ominus}}$.}
\end{Prop}

\subsubsection{Comparaison des résultats obtenus}

Remarquons tout d'abord que le rapprochement des Propositions \ref{P6} et \ref{P6'} donne bien l'identité:\smallskip

\centerline{$q_\Delta(G,\Cl^{\,\si{\ominus}})+q_\Delta(G,\C^{\,\si{\ominus}})=0$.}\smallskip

Cependant le calcul du premier quotient $q_\Delta(G,\Cl^{\,\si{\ominus}})$ vaut en toute généralité, tandis que celui du second $q_\Delta(G,\C^{\,\si{\ominus}})$ présuppose la trivialité de l'invariant mu d'Iwasawa, qui garantit la finitude du groupe $H^1(G,\C_{L'_\%}^{\,\si{\ominus}})$.
D'un autre côté, dans le calcul du caractère du groupe $H^2(G,\C_{L'_\%}^{\,\si{\ominus}})$, nous avons utilisé en outre un résultat sur les diviseurs imaginaires au-dessus de $\ell$ dans une $\Zl$-extension cyclotomique, qui est une conséquence facile des conjectures de Leopoldt et de Gross-Kuz'min (cf. \cite{J18}, Th. IV.2.19), mais dont il peut être intéressant de donner une démonstration autonome:

\begin{Prop}\label{P7}
Soient $L'_\%=\bigcup_{n\in\NN}L'_n$ la $\Zl$-extension cyclotomique associée à un corps de nombres à conjugaison complexe $L'$ contenant les racines $\ell$-ièmes de l'unité, et $n_{\si{0}}$ le rang à partir duquel les places au-dessus de $\ell$ se ramifient totalement dans la tour. Alors:
\begin{itemize}\smallskip

\item[(i)] Pour tout $n \geqslant n_{\si{0}}$ les $\ell$-unités imaginaires du corps $L'_n$ sont les produits des $\ell$-unités de $L'_{n_{\si{0}}}$ par les racines de l'unité contenues dans $L'_n$; ce qui s'écrit:\smallskip

\centerline{$E'_{L'_n}{\!\!\!\!}^{\si{\ominus}}=\mu_{L'_n}E'_{L'_{n_{\si{0}}}}{\!\!\!\!\!\!}^{\si{\ominus}}$.}\smallskip

\item[(ii)] Les diviseurs imaginaires principaux de $L'_n$ construits sur les places au-dessus de $\ell$ sont les étendus à $L'_n$ des diviseurs imaginaires principaux de $L'_{n_{\si{0}}}$ construits sur ces mêmes places.\smallskip

\item[(iii)] La composante imaginaire  $\,\Cl_{L'_\%}^{[\ell]\si{\,\ominus}}$ de la limite inductive des $\ell$-sous-groupes de classes des corps $L'_n$ construit sur les classes des places au-dessus de $\ell$ est un $\Zl$-module divisible qui a pour codimension le nombre $d$ de places de $L'_\%$ au-dessus de $\ell$ décomposées par la conjugaison complexe.\smallskip

\item[(iv)] La composante imaginaire $\,\C_{L'_\%}^{[\ell]\,\si{\ominus}}$ de la limite projective est un $\Zl$-module libre de rang $d$.
\end{itemize}
\end{Prop}

\Preuve 
Désignons par $D_{L'_\%}^{[\ell]\,\si{\ominus}}$ le groupe des diviseurs imaginaires de $L'_\%$ construits sur les places au-dessus de $\ell$. D'après la suite exacte courte canonique\smallskip

\centerline{$1 \to E'{\!}^{\;\si{\ominus}}_{\!L'_\%}/\mu^\ph_{L'_\%} \to D_{L'_\%}^{[\ell]\,\si{\ominus}} \to \Cl_{L'_\%}^{[\ell]\,\si{\ominus}} \to 1$}\smallskip

\noindent il s'agit de vérifier que le quotient $E'{\!}^{\;\si{\ominus}}_{\!L'_n}/\mu^\ph_{L'_n}$ à gauche est constant pour $n\geqslant n_{\si{0}}$. Or, cela résulte d'un argument classique de théorie de Kummer (cf. \cite{J18}, Prop. I.2.2): l'égalité des rangs montre que  $E'{\!}^{\;\si{\ominus}}_{\!L'_n}\mu^\ph_{L'_n}$ est d'indice fini dans  $E'{\!}^{\;\si{\ominus}}_{\!L'_{n+\si{1}}}/\mu^\ph_{L'_{n+\si{1}}}$ pour $n\geqslant n_{\si{0}}$. Cela étant, si $x$ est un élément de $E'{\!}^{\;\si{\ominus}}_{\!L'_n}$ dont la classe modulo $\mu^\ph_{L'_n}$ est une puissance $\ell$-ième dans $E'{\!}^{\;\si{\ominus}}_{\!L'_{n+\si{1}}}/\mu^\ph_{L'_{n+\si{1}}}$, nous avons $x=\zeta y^\ell$ avec $\zeta\in\mu^\ph_{L'_{n+\si{1}}}$ et $y\in E'{\!}^{\;\si{\ominus}}_{\!L'_{n+\si{1}}}$, puis, en prenant la norme $N=N_{\si{L'_{n+\si{1}}/L'_n}}$, l'égalité $x^\ell=\zeta^\ell N(y)^\ell$, i.e. $x/\zeta N(y) \in \mmu_\ell \subset L'_n$; donc $\zeta \in L'_n$, et finalement: $\zeta \in \mu_{L'_n}^\ph=\mu_{L'_{n+\si{1}}}^\ell$. Écrivant donc $\zeta=\xi^\ell$, et $z=\xi y$, nous obtenons $x=z^\ell$, donc $L'_n[\sqrt[\ell]{x}]=L'_n[z]\subset L'_{n+\si{1}}$; et la théorie de Kummer nous dit que $x$ est le produit dans $L'_n$ d'une puissance $\ell$-ième et d'une racine de l'unité. Et  $E'{\!}^{\;\si{\ominus}}_{\!L'_n}/\mu^\ph_{L'_n}$ coïncide avec $E'{\!}^{\;\si{\ominus}}_{\!L'_{n+\si{1}}}/\mu^\ph_{L'_{n+\si{1}}}$ pour $n\geqslant n_{\si{0}}$. La Proposition en découle.
\medskip

Revenons maintenant sur le quotient $q_\Delta(G,\,\Cl^{\;\si{\ominus}}_{\!L'_\%})$. Dans le diagramme commutatif
\begin{displaymath}
\xymatrix{
1 \ar@{->}[r] & H^1(G,\Cl^{\;\si{\ominus}}_{L'_\%}) \ar@{->}[r] \ar@{->}[d]^\simeq & H^2(G,P^{\;\si{\ominus}}_{L'_\%}) \ar@{->}[r] \ar@{->}[d]^\simeq & H^2(G,D^{\;\si{\ominus}}_{L'_\%}) \ar@{->}[r] \ar@{->}[d]^\simeq & H^2(G,\Cl^{\;\si{\ominus}}_{L'_\%}) \ar@{->}[r] \ar@{->}[d]^\simeq & 1\\
1 \ar@{->}[r] & \C\!ap^{\;\si{\ominus}}_{L'_\%/K'_\%} \ar@{->}[r] & H^1(G,\mmu_\ell^\ph) \ar@{->}[r]  & \Ram '{}^{\si{\ominus}}_{\!\!L'_\% / K'_\%}  \ar@{->}[r] &  \Cl_{L'_\%}^{\,G\,\si{\ominus}}  / \Cl^{\;\et\,\si{\ominus}}_{L'_\%} \ar@{->}[r] & 1 
}
\end{displaymath}
l'isomorphisme à gauche, induit par la norme arithmétique $N_{\si{L'_\%/K'_\%}}$, identifie le groupe $H^1(G,\Cl^{\;\si{\ominus}}_{\!L'_\%})$ à la capitulation imaginaire $\C\!ap^{\;\si{\ominus}}_{L'_\%/K'_\%}$ dans $L'_\%/K'_\%$. Comme le groupe $H^1(G,\mmu_\ell^\ph)\simeq\mmu_\ell$ est un $\Fl[\Delta]$-module de caractère $\omega$, deux cas peuvent se produire:
\begin{itemize}
\item ou bien $\,\C\!ap^{\;\si{\ominus}}_{L'_\%/K'_\%}=1$ et $|H^2(G,\Cl^{\;\si{\ominus}}_{\!L'_\%})| = |\Ram '{}^{\si{\ominus}}_{\!\!L'_\% / K'_\%}|/\ell$;
\item ou bien $\,\C\!ap^{\;\si{\ominus}}_{L'_\%/K'_\%} \simeq H^1(G,\mmu_\ell^\ph)$ et $H^2(G,\Cl^{\;\si{\ominus}}_{\!L'_\%})\simeq \Ram '{}^{\si{\ominus}}_{\!\!L'_\% / K'_\%}$.
\end{itemize}
Le premier cas nécessite que le groupe $\Ram '{}^{\si{\ominus}}_{\!\!L'_\% / K'_\%}$ soit non-trivial, i.e. qu'il existe au moins une place $\p_\%^{\si{+}}$ étrangère à $\ell$, décomposée dans $LK\%/K_\%^+$ qui se ramifie dans $L^+_\%/K^+_\%$. Le second cas implique, en revanche, que l'extension $L^+_\%/K^+8_\%$ soit $\ell$-ramifiée: en effet, si $\a$ est un diviseur de $K'_\%$ représentant une classe de $\,\Cl^{\;\si{\ominus}}_{K'_\%}$ qui capitule dans $\,\Cl^{\;\si{\ominus}}_{L'_\%}$, son étendu à $L'_\%$ est un diviseur principal $(\alpha)$ engendré par un élément imaginaire $\alpha$ qui vérifie $(\alpha)^{\sigma-\si{1}}=1$, i.e. $\alpha^{\sigma-\si{1}}\in \mmu_\ell$. Posant $a=\alpha^\ell$, nous obtenons donc $L'_\%=K'_\%[\sqrt[\ell]{a}]$, avec $(a)=\a^\ell$, de sorte que l'extension $L'_\%/K'_\%$ est $\ell$-ramifiée.
L'ensemble des résultats précédents peut donc se résumer comme suit:

\begin{Th}\label{T2}
Soient $L_\%/K_\%$ une $\ell$-extension cyclique élémentaire de corps surcirculaires à conjugaison complexe, $G$ son groupe de Galois, $L^+_\%/K^+_\%$ la sous-extension réelle maximale et $L'_\%/K'_\%=L_\%[\zeta_\ell]/K_\%[\zeta_\ell]$ l'extension obtenue en adjoignant les racines $\ell$-ièmes de l'unité, $\Delta$ enfin le groupe abélien $\Gal(L'_\%/L^+_\%)\simeq\Gal(K'_\%/K^+_\%)$.
Alors, sous réserve de nullité de l'invariant mu, les caractères de $\Delta$ attachés aux groupes de cohomologie relatifs à l'action de $G$ sur les composantes imaginaires du $\ell$-groupe des classes de diviseurs $\,\Cl_{L'_\%}$ et du groupe de Galois $\,\C_{L'_\%}$ de la pro-$\ell$-extension abélienne non-ramifiée maximale de $L'_\%$ sont donnés par les formules:\smallskip

$\begin{aligned}
&(i) \quad h^1_\Delta(G,\Cl^{\,\si{\ominus}}_{L'_\%}) =\delta'\omega \phantom{AAAAAAAAAAAA}  &(i)'& \quad h^2_\Delta(G,\C_{L'_\%}^{\,\si{\ominus}})= \sum_{\p_\%^{\si{+}}\mid\ell}\chi_{\p_\%^{\si{+}}}^{\si{\ominus}}+(1-\delta)\omega \\
&(ii) \quad h^2_\Delta(G,\Cl^{\,\si{\ominus}}_{L'_\%})=\sum_{\p_\%^{\si{+}}\mid\ell}\chi_{\p_\%^{\si{+}}}^{\si{\ominus}}+(1-\delta')\omega  &(ii)'& \quad h^1_\Delta(G,\C^{\,\si{\ominus}}_{L'_\%})=\sum_{\p_\%^{\si{+}}}\chi_{\p_\%^{\si{+}}}^{\si{\ominus}}-\delta\omega
\end{aligned}$

\noindent Dans celles-ci, $\omega$ est le caractère du module de Tate $\Tl=\varprojlim \mmu_{\ell^{\si{n}}}=\zeta^\Zl$; l'indice $\delta$ vaut 0 ou 1 suivant que $\zeta$ est norme ou non dans l'extension $L'_\%/K'_\%$ ; l'indice $\delta'$ vaut 0 ou 1 suivant qu'il existe ou non une place de $K^+_\%$ décomposée dans $K_\%/K_\%^+$ qui se ramifie modérément dans $L^+_\%/K_\%^+$; et $\chi_{\p_\%^{\si{+}}}^{\si{\ominus}}$ désigne la composante imaginaire de l'induit à $\Delta$ du caractère unité du sous-groupe de décomposition de $\p_\%^{\si{+}}$ dans l'extension abélienne $K'_\%/K^+_\%$.
\end{Th}

\begin{Sco}\label{S}
Sous les hypothèses du Théorème \ref{T2}, l'identité $h^i_\Delta(G,\C^{\,\si{\ominus}}_{L'_\%})=h^{i+{\si{1}}}_\Delta(G,\Cl^{\,\si{\ominus}}_{L'_\%})$ a lieu exactement dans les deux cas suivants:
\begin{itemize}
\item[(i)] Dans le cas modéré, i.e. lorsque toutes les places de $K^+_\%$, décomposées dans l'extension quadratique $K_\%/K_\%^+$ et ramifiées dans l'extension cyclique $L_\%/K_\%$, sont étrangères à $\ell$.
\item[(ii)] Dans un cas sauvage exceptionnel, lorsqu'il existe exactement une place $\p_\%^{\si{+}}$ de $K^+_\%$ décomposée dans $K_\%/K_\%^+$, qui se ramifie dans l'extension $L_\%/K_\%$: si $\p_\%^{\si{+}}$ est au-dessus de $\ell$, sous réserve qu'il y ait exactement deux places $\p_\%'$ et $\bar\p_\%'$ de $K'_\%$ au-dessus de $\p_\%^{\si{+}}$ et que $\zeta$ ne soit norme locale dans $L'_\%/K'_\%$ en aucune d'elles.
\end{itemize}
\end{Sco}

\Preuve Écrivons $\Sigma$ pour $ \sum_{\p_\%^{\si{+}}\mid\ell}\chi_{\p_\%^{\si{+}}}^{\si{\ominus}}$. D'après le théorème \ref{T2}, les cohomologies de $\,\C^{\,\si{\ominus}}_{L'_\%}$ et de $\,\Cl^{\,\si{\ominus}}_{L'_\%}$ sont duales lorsqu'on a:
$\Sigma + (1-\delta-\delta')\omega=0$, i.e. 
$\Sigma=0$ et $\delta+\delta'=1$ ou bien $\Sigma=\omega$ et $\delta=\delta'=1$.\smallskip

Dans le premier cas ($\Sigma=0$), nous disons que la ramification imaginaire est modérée. Suivant qu'il existe ou non une place $\p_\%^{\si{+}}$ de $K^+_\%$, décomposée dans $K_\%/K_\%^+$ qui se ramifie modérément dans $L_\%^+/K_\%^+$, nous avons alors $\delta'=1$ et $\delta=0$ ou $\delta'=0$ et $\delta=1$\footnote{Conformément au principe de Hasse, la condition $\delta=0$ se lit localement. Elle est automatiquement vérifiée en l'absence de place ramifiée. En revanche, en présence d'une place modérément ramifiée, le corps de classes local montre que les racines primitives de l'unité ne peuvent être normes, ce qui donne $\delta=1$.}\; donc toujours: $\delta+\delta'=1$.

\subsection{Énoncé des résultats}

Revenons maintenant au problème posé dans l'introduction: désignons par $L_\%/K_\%$ une $\ell$-extension de corps surcirculaires à conjugaison complexe; notons $G$ son groupe de Galois, $L^+_\%/K^+_\%$ la sous-extension réelle maximale et $L'_\%/K'_\%=L_\%[\zeta_\ell]/K_\%[\zeta_\ell]$ l'extension obtenue en adjoignant les racines $\ell$-ièmes de l'unité, $\Delta$ enfin le groupe abélien $\Gal(L'_\%/L^+_\%)\simeq\Gal(K'_\%/K^+_\%)$.\smallskip

Pour chaque sous-extension $M'_\%/K'_\%$ de $L'_\%/K'_\%$, nous savons par la Proposition \ref{P3}, sous réserve de nullité de l'invariant mu, que la composante imaginaire $\,\C_{M'_\%}^{\,\si{\ominus}}$ du groupe de Galois de la pro-$\ell$-extension abélienne non-ramifiée maximale de $M'_\%$ est un $\Zl[\Delta]$-module projectif dont nous avons convenu d'écrire $\lambda_{M'_\%}^{\,\si{\ominus}}$ le caractère (de sorte que c'est son degré $\deg\lambda_{M'_\%}^{\,\si{\ominus}}$ qui désigne l'invariant lambda d'Iwasawa de $\,\C_{M'_\%}^{\,\si{\ominus}}$). 
Si, en outre, $M'_\%/K'_\%$ est galoisienne, le groupe $\,\C_{M'_\%}^{\,\si{\ominus}}$ est muni d'une action naturelle de $G$ et nous écrivons $\chi_{M'_\%}^{\,\si{\ominus}}$ le caractère de $\,\C_{M'_\%}^{\,\si{\ominus}}$ regardé comme $\Delta \times G$-module.\smallskip

Notre propos est d'exprimer $\chi_{L'_\%}^{\,\si{\ominus}}$ à partir de $\chi_{K'_\%}^{\,\si{\ominus}}$.\smallskip

Restreignons-nous un instant au cas élémentaire où le groupe $G$ et cyclique d'ordre premier $\ell$. Dans ce cas, la Proposition \ref{P5} nous donne l'identité:\smallskip

\centerline{$\chi_{L'_\%}^{\,\si{\ominus}} = \chi_{K'_\%}^{\,\si{\ominus}}\Reg_G - q(G,\C_{L'_\%}^{\,\si{\ominus}})\Aug_G$}\smallskip

\noindent où $\Reg_G$ est le caractère régulier de $G$ et $\Aug_G$ le caractère d'augmentation. Quant au quotient de Herbrand dimensionnel $ q(G,\C_{L'_\%}^{\,\si{\ominus}})$, il est donné par la Proposition \ref{P6'}: il est égal à la différence\smallskip

\centerline{$\omega - \sum_{\p_\%^{\si{+}}\nmid\ell}\chi_{\p_\%^{\si{+}}}^{\si{\ominus}}$}\smallskip

\noindent où la sommation porte sur les places $\p_\%^{\si{+}}$de $K^+_\%$ étrangères à $\ell$ qui se ramifient dans $L_\%^+/K_\%^+$.\smallskip

Nous pouvons donc réécrire l'identité précédente sous la forme:\smallskip

\centerline{$\chi_{L'_\%}^{\,\si{\ominus}}-\omega = (\chi_{K'_\%}^{\,\si{\ominus}}-\omega)\Reg_G + \sum_{\p_\%^{\si{+}}\nmid\ell}\,\chi_{\p_\%^{\si{+}}}^{\si{\ominus}}\,\rho_{\p_\%^{\si{+}}}^\ph$}\smallskip

\noindent en convenant de désigner par $\rho_{\p_\%^{\si{+}}}^\ph$ l'induit à $G\simeq\Gal(L_\%/K_\%^+)$ du caractère d'augmentation du sous-groupe de décomposition $G_{\p_\%^{\si{+}}}^\ph$ de $\p_\%^{\si{+}}$ dans $L_\%^+/K_\%^+$. 
Dans cette dernière formulation, la sommation à droite porte indifféremment sur les seules places  étrangères à $\ell$ ramifiées dans $L_\%^+/K_\%^+$ ou sur l'ensemble des places de $K_\%^+$ étrangères à $\ell$, puisque le caractère  $,\rho_{\p_\%^{\si{+}}}^\ph$, égal à $\Aug_G$ lorsque $\p_\%^{\si{+}}$ se ramifie dans $L_\%^+/K_\%^+$, est nul dans le cas contraire faute d'inertie aux places étrangères à $\ell$.\smallskip

Maintenant, comme expliqué au début de cette section, la formule ainsi obtenue, établie dans le cas où $G$ est d'ordre $\ell$, vaut en fait pour un $\ell$-groupe arbitraire. Nous pouvons donc énoncer:

\begin{Th}\label{T3}
Soient $L_\%/K_\%$ une $\ell$-extension cyclique élémentaire de corps surcirculaires à conjugaison complexe, $G$ son groupe de Galois, $L^+_\%/K^+_\%$ la sous-extension réelle maximale et $L'_\%/K'_\%=L_\%[\zeta_\ell]/K_\%[\zeta_\ell]$ l'extension obtenue en adjoignant les racines $\ell$-ièmes de l'unité, $\Delta$ enfin le groupe abélien $\Gal(L'_\%/L^+_\%)\simeq\Gal(K'_\%/K^+_\%)$.
Alors, sous réserve de nullité de l'invariant mu d'Iwasawa dans $L'_\%$, les composantes imaginaires respectives $\,\C_{L'_\%}^{\,\si{\ominus}}$ et $\,\C_{K'_\%}^{\,\si{\ominus}}$ des groupes de Galois des pro-$\ell$-extensions abéliennes non-ramifiées maximales de $L'_\%$ et $K'_\%$ sont des $\Zl[\Delta]$-modules projectifs et leurs caractères  $\chi_{L'_\%}^{\,\si{\ominus}}$ et $\chi_{K'_\%}^{\,\si{\ominus}}$ comme $\Delta\times G$-modules sont liés par l'identité:\smallskip

\centerline{$\chi_{L'_\%}^{\,\si{\ominus}} -\omega = (\chi_{K'_\%}^{\,\si{\ominus}} -\omega)\Reg_G + \sum_{\p_\%^{\si{+}}\nmid\ell}\chi_{\p_\%^{\si{+}}}\,\rho_{\p_\%^{\si{+}}}^\ph$.}\smallskip

\noindent Dans celle-ci, $\omega$ est le caractère cyclotomique; $\Reg_G$ est le caractère régulier de $G$; la sommation porte sur les places $\p_\%^+$ de $K^+_\%$ qui sont étrangères à $\ell$; le caractère $\chi_{\p_\%^{\si{+}}}$ est l'induit à $\Delta$ du caractère unité du sous-groupe de décomposition $\Delta_{\chi_{\p_\%^{\si{+}}}}$ de $\p_\%^{\si{+}}$ dans $K'_\%/K^+_\%$ et $\rho_{\p_\%^{\si{+}}}$ est l'induit à $G$ du caractère d'augmentation du sous-groupe de décomposition de $\p_\%^{\si{+}}$ dans $L^+_\%/K^+_\%$.
\end{Th}

\begin{Cor}[Formule d'Iwasawa]\label{C1}
Notons $\Delta'=\Gal(L'_\%/L_\%)\simeq\Gal(K'_\%/K_\%)$. Avec les notations du Théorème \ref{T3}, les caractères respectifs des groupes $\,\C_{L'_\%}^{\,\si{\ominus}}$ et $\,\C_{K'_\%}^{\,\si{\ominus}}$ pour leur structure de $G$-module sont liés par l'identité:\smallskip

\centerline{$\chi_{L_\%}^{\,\si{\ominus}} -\langle\omega,1_{\Delta'}^\Delta\rangle 1_G = (\chi_{K'_\%}^{\,\si{\ominus}} -\langle\omega,1_{\Delta'}^\Delta\rangle 1_G)\Reg_G + ( \sum_{\p_\%^{\si{+}}\nmid\ell}\langle\chi_{\p_\%^{\si{+}}},1_{\Delta'}^\Delta\rangle)\,\rho_{\p_\%}^\ph$.}\smallskip

Dans celle-ci, $1_{\Delta'}^\Delta$ est l'induit à $\Delta$ du caractère unité de $\Delta'$; la quantité $\langle\omega,1_{\Delta'}^\Delta\rangle$ vaut 1 si $K_\%$ contient le racines $\ell$-ièmes de l'unité et 0 sinon; et la quantité $\langle\chi_{\p_\%^{\si{+}}},1_{\Delta'}^\Delta\rangle$ vaut 1 ou 0 selon que la place ${\p_\%^{\si{+}}}$ est décomposée ou non par la conjugaison complexe.
\end{Cor}

\Preuve
Les invariants attachés aux corps $L_\%$ et $K_\%$ s'obtiennent à partir de ceux respectivement attachés à $L'_\%$ et $K'_\%$ par projection au moyen de l'idempotent associé au caractère $1_{\Delta'}^\Delta$.

\begin{Cor}[Formule de Kida]\label{A1}
Sous les hypothèses du Théorème \ref{T3}, les invariants lambda d'Iwasawa des groupes $\,\C_{L'_\%}^{\,\si{\ominus}}$ et $\,\C_{K'_\%}^{\,\si{\ominus}}$ sont liés par l'identité:\smallskip

\centerline{$\lambda_{L_\%}^{\,\si{\ominus}} -\langle\omega,1_{\Delta'}^\Delta\rangle = (\lambda_{K_\%}^{\,\si{\ominus}} -\langle\omega,1_{\Delta'}^\Delta\rangle)[L_\%:K_\%] + \sum_{\si{\mathfrak P}_\%^{\si{+}}\nmid\ell}(e_{\si{\mathfrak P}_\%^{\si{+}}}-1)$,}\smallskip

\noindent où la sommation porte sur les places $\si{\mathfrak P}_\%^{\si{+}}$ de $L^+_\%$ modérément ramifiées dans $L^+_\%/K^+_\%$ et décomposées dans $L_\%/L^+_\%$; et $e_{\si{\mathfrak P}_\%^{\si{+}}}$ est l'indice de ramification de $\si{\mathfrak P}_\%^{\si{+}}$ dans $L^+_\%/K^+_\%$.
\end{Cor}

\begin{Cor}[Formule de Wingberg]\label{A2}
Sous les mêmes hypothèses, les invariants lambda des $\ell$-groupes de classes infinitésimales des corps totalement réels $L_\%+$ et $K_\%^+$ sont lié, sous la conjecture de Leopoldt, par l'identité:\smallskip

\centerline{$\lambda_{\{\ell\}\,L^+_\%} -\langle\omega,1_{\Delta'}^\Delta\rangle = (\lambda_{\{\ell\}\,K^+_\%} -\langle\omega,1_{\Delta'}^\Delta\rangle)[L_\%:K_\%] + \sum_{\si{\mathfrak P}_\%^{\si{+}}}(e_{\si{\mathfrak P}_\%^{\si{+}}}-1)$.}
\end{Cor}

\Preuve
Prenant les degrés des divers caractères dans la formule d'Iwasawa, on obtient la formule de Kida. D'après la Proposition \ref{P3}, la même identité vaut sous la conjecture de Leopoldt pour les invariants $\lambda_{\si{\{\ell\}}}$ en vertu de l'égalité $\lambda_{\si{\{\ell\}}}^{\,\si{\oplus}}=\lambda^{\si{\ominus}}_\ph$; ce qui redonne la formule de Wingberg.\medskip

Revenons enfin sur le Théorème \ref{T3} et remplaçons $\,\C_{L'_\%}$ par le groupe de Galois $\,\wC_{L'_\%}$ de la pro$-\ell$-extension abélienne non-ramifiée $\ell$-décomposée du corps $L'_\%$; regardons-le comme $\Delta\times G$-module et notons $\wi\chi_{L'_\%}^{\,\si{\ominus}}$ le caractère de sa composante imaginaire. Par le corps de classes, $\,\wC_{L'_\%}$ s'identifie à la limite projective $\varprojlim \Cl'_{L'_n}$ des $\ell$-groupes de $\ell$-classes des corps $L'_n$, autrement dit au quotient $\,\C_{L'_\%}/\C_{L'_\%}^{\,\si{[\ell]}}$ de $\,\C_{L'_\%}$ par le sous-groupe $\,\C_{L'_\%}^{\,\si{[\ell]}}=\varprojlim \Cl_{L'_n}^{\,\si{[\ell]}}$ construit sur les classes des diviseurs au-dessus de $\ell$.
D'après la Proposition \ref{P7}, les caractères $\chi^{\,\si{\ominus}}$ et $\wi\chi^{\,\si{\ominus}}$ sont donc liés par les identités:\smallskip

\centerline{$\chi_{K'_\%}^{\,\si{\ominus}}=\wi\chi_{K'_\%}^{\,\si{\ominus}} + \sum_{\p_\%^{\si{+}}\mid\ell} \chi_{\p_\%^{\si{+}}}^{\,\si{\ominus}}\,1^\ph_G$ 
\qquad et \qquad $\chi_{L'_\%}^{\,\si{\ominus}}=\wi\chi_{L'_\%}^{\,\si{\ominus}} + \sum_{\p_\%^{\si{+}}\mid\ell} \chi_{\p_\%^{\si{+}}}^{\,\si{\ominus}}\,1_{G_{\si{\mathfrak P}_\%^{\si{+}}}}^G$,}

\noindent où $1_{G_{\si{\mathfrak P}_\%^{\si{+}}}}^G$ désigne l'induit à $G$ du caractère unité du sous-groupe de décomposition $G_{\si{\mathfrak P}_\%^{\si{+}}}$ dans $L_\%^+/K_\%^+$ de l'une quelconque $\si{\mathfrak P}_\%^{\si{+}}$ des places de $L_\%^+$ au-dessus de $\p_\%^{\si{+}}$.
\smallskip

Par différence, nous avons donc:\smallskip

\centerline{$\chi_{L'_\%}^{\,\si{\ominus}} -  \chi_{K'_\%}^{\,\si{\ominus}}\Reg^\ph_G = \wi\chi_{L'_\%}^{\,\si{\ominus}} -  \wi\chi_{K'_\%}^{\,\si{\ominus}}\Reg^\ph_G + \sum_{\p_\%^{\si{+}}\mid\ell} \chi_{\p_\%^{\si{+}}}^{\,\si{\ominus}}\,\rho^\ph_{\p_\%^{\si{+}}}$,}\smallskip

\noindent où $\rho^\ph_{\p_\%^{\si{+}}}=\Aug_{G_{\si{\mathfrak P}_\%^{\si{+}}}}^G$ est l'induit à $G$ du caractère d'augmentation de  $G_{\si{\mathfrak P}_\%^{\si{+}}}$. Il suit:

\begin{Th}\label{T3'}
Sous les hypothèses et avec les notations du Théorème \ref{T3}, les composantes imaginaires respectives $\,\wC_{L'_\%}^{\,\si{\ominus}}$ et $\,\wC_{K'_\%}^{\,\si{\ominus}}$ des groupes de Galois des pro-$\ell$-extensions abéliennes non-ramifiées $\ell$-décomposées maximales de $L'_\%$ et $K'_\%$ sont des $\Zl[\Delta]$-modules projectifs et leurs caractères  $\wi\chi_{L^+_\%}^\ph$ et $\wi\chi_{K^+_\%}^\ph$ comme $\Delta\times G$-modules sont liés par l'identité:\smallskip

\centerline{$\wi\chi_{L^+_\%}^{\,\si{\ominus}} -\omega = (\wi\chi_{K^+_\%}^{\,\si{\ominus}} -\omega)\Reg_G + ( \sum_{\p_\%^{\si{+}}} \chi_{\p_\%^{\si{+}}}^{\,\si{\ominus}})\,\rho^\ph_{\p_\%^{\si{+}}}$.}\smallskip

\noindent où la sommation porte ici sur l'ensemble des places $\p_\%^{\si{+}}$ du corps $K_\%^+$.
\end{Th}

\begin{Cor}\label{C2}
Toujours avec les mêmes notations, les caractères respectifs des groupes $\,\wC_{L_\%}^{\,\si{\ominus}}$ et $\,\wC_{K_\%}^{\,\si{\ominus}}$ regardés comme $G$-modules sont liés par l'identité:\smallskip

\centerline{$\wi\chi_{L_\%}^{\,\si{\ominus}} -\langle\omega,1_{\Delta'}^\Delta\rangle 1_G = (\wi\chi_{K_\%}^{\,\si{\ominus}} -\langle\omega,1_{\Delta'}^\Delta\rangle 1_G)\Reg_G + ( \sum_{\p_\%^{\si{+}}}\langle\chi_{\p_\%^{\si{+}}}^{\,\si{\ominus}},1_{\Delta'}^\Delta\rangle)\,\rho_{\p_\%}^{\si{+}}$.}
\end{Cor}

\begin{Cor}[Formule de Kuz'min]\label{A3}
Sous les mêmes hypothèses, les invariants lambda des composantes imaginaires des $\ell$-groupes $\,\wC_{L_\%}$ et $\,\wC_{K_\%}$ sont liés par l'identité:\smallskip

\centerline{$\wi\lambda_{L_\%}^{\,\si{\ominus}} -\langle\omega,1_{\Delta'}^\Delta\rangle = (\wi\lambda_{K_\%}^{\,\si{\ominus}} -\langle\omega,1_{\Delta'}^\Delta\rangle)[L_\%:K_\%] + \sum_{\si{\mathfrak P}_\%^{\si{+}}}(d_{\si{\mathfrak P}_\%^{\si{+}}}-1)$,}\smallskip

\noindent où la sommation porte sur les places $\si{\mathfrak P}_\%^{\si{+}}$ de $L^+_\%$ modérément ramifiées dans $L^+_\%/K^+_\%$ et décomposées dans $L_\%/L^+_\%$; et $d_{\si{\mathfrak P}_\%^{\si{+}}}=e_{\si{\mathfrak P}_\%^{\si{+}}}$ est l'ordre du sous-groupe de décomposition $G_{\si{\mathfrak P}_\%^{\si{+}}}$ de $\si{\mathfrak P}_\%^{\si{+}}$ dans $L^+_\%/K^+_\%$.
\end{Cor}

\newpage
\def\refname{\normalsize{\sc  Références}}
\addcontentsline{toc}{section}{Bibliographie}
{\small

\def\refname
}

\bigskip
{\small
\noindent{\sc Adresse:}
Institut de Mathématiques de Bordeaux, \par
351 Cours de la Libération,
F-33405 Talence cedex

\noindent{\sc Courriel:}
{\tt jean-francois.jaulent@math.u-bordeaux.fr}
 
 \noindent{\sc url:} \url{https://www.math.u-bordeaux.fr/~jjaulent/}
}

\end{document}